\documentclass{article}

\usepackage[preprint]{neurips_2024}

\usepackage[utf8]{inputenc} 
\usepackage[T1]{fontenc}    
\usepackage{hyperref}       
\usepackage{url}            
\usepackage{booktabs}       
\usepackage{amsfonts}       
\usepackage{nicefrac}       
\usepackage{microtype}      
\usepackage{xcolor}         

\usepackage{style_bo}
\usepackage{amsthm}
\newtheorem{theorem}{Theorem}
\newtheorem{remark}{Remark}
\newtheorem{proposition}{Proposition}
\newtheorem{definition}{Definition}
\newtheorem{lemma}{Lemma}

\newtheorem{example}{Example}
\newtheorem{assumption}{Assumption}

\DeclareMathOperator*{\argmin}{argmin}

\title{Conformal Inverse Optimization}

\author{%
  Bo Lin \\
  University of Toronto \\
  \texttt{blin@mie.utoronto.ca} \\
  \And
  Erick Delage \\
  HEC Montréal \\
  \texttt{erick.delage@hec.ca} \\
  \And
  Timothy C. Y. Chan \\
  University of Toronto \\
  \texttt{tcychan@mie.utoronto.ca}
}

\begin{document}

\maketitle

\begin{abstract}
Inverse optimization has been increasingly used to estimate unknown parameters in an optimization model based on decision data. We show that such a point estimation is insufficient in a prescriptive setting where the estimated parameters are used to prescribe new decisions. The prescribed decisions may be low-quality and misaligned with human intuition and thus are unlikely to be adopted. To tackle this challenge, we propose conformal inverse optimization, which seeks to learn an uncertainty set for the unknown parameters and then solve a robust optimization model to prescribe new decisions. Under mild assumptions, we show that our method enjoys provable guarantees on solution quality, as evaluated using both the ground-truth parameters and the decision maker's perception of the unknown parameters. Our method demonstrates strong empirical performance compared to classic inverse optimization.
\end{abstract}

\section{Introduction}

Inverse optimization (IO) is a supervised learning approach that fits parameters in an optimization model to decision data. The fitted optimization model can then be used to prescribe future decisions. Such problems naturally arise in AI applications where human preferences are not explicitly given, and instead need to be inferred from historical decisions. For this pipeline to succeed in practice, the prescribed decision should not only be of high-quality but also align with human intuition (i.e., perceived to be of high-quality). The latter encourages algorithm adoption \citep{chen2023understanding, donahue2023two}, which is critical in many real-world settings \citep{liu2023algorithm, sun2022predicting}. 

As an example, rideshare platforms, e.g., Uber and Lyft, provide a shortest-path to the driver at the start of a trip based on real-time traffic data \citep{nguyen2015uber}. The driver then relies on her perception of the road network formed through past experience to evaluate the path. The driver may deviate from the suggested path if it is perceived to be low-quality. Although seasoned drivers are often capable of identifying a better path due to their tacit knowledge of the road network \citep{merchan20222021}, such deviations impose operational challenges as it may cause rider safety concerns and affect downstream decisions such as arrival time estimation, trip pricing, and rider-driver matching \citep{hu2022deepeta}. Therefore, the platform may be interested in leveraging historical paths taken by drivers to suggest high-quality paths for future trips, as evaluated using both the travel time and the driver's perception. 

In this paper, we first show that the classic IO pipeline may generate decisions that are low-quality and misaligned with human intuition. We next propose conformal IO, which first learns an uncertainty set from decision data and then solves a robust optimization model with the learned uncertainty set to prescribe decisions. Finally, we show that the proposed approach has provable guarantees on the actual and perceived solution quality. Our contributions are as follows.

\textbf{New framework}. We propose a new IO pipeline that integrates i) a novel method to learn uncertainty sets from decision data and ii) a robust optimization model for decision recommendation.

\textbf{Theoretical guarantees}. We prove that, with high probability, the learned uncertainty set contains parameters that make future observed decisions optimal. This coverage guarantee leads to bounds on the optimality gap of the decisions from conformal IO, as evaluated using the ground-truth parameters and the decision maker's (DM's) perceived parameters.

\textbf{Performance}. Through experiments, we demonstrate strong empirical performance of conformal IO compared to classic IO and provide insights into modeling choices.

\section{Literature Review}

\textbf{Inverse optimization}. IO is a method to estimate unknown parameters in an optimization problem based on decision data \citep{ahuja2001inverse, chan2014generalized, chan2020inverse}. Early IO papers focus on deterministic settings where the observed decisions are assumed to be optimal to the specified optimization model. Recently, IO has been extended to stochastic settings where the observed decisions are subject to errors and bounded rationality. Progress has been made to provide estimators that are consistent \citep{aswani2018inverse, dong2018generalized}, tractable \citep{chan2019inverse,  tan2020learn}, and robust to data corruption \citep{mohajerin2018dataIO}. Our paper is in the stochastic stream. Unlike existing methods that provide a point estimation of the unknown parameters, we learn an uncertainty set that can be used in a robust optimization model.

\textbf{Data-driven uncertainty set construction.} Recently, data have become a critical ingredient to design the structure \citep{delage2010distributionally, mohajerin2018dataIO, gao2023distributionally} and calibrate the size \citep{chenreddy2022data, sun2023predict} of an uncertainty set. Our paper is related to the work of \citet{sun2023predict} who first use an ML model to predict the unknown parameters and then calibrate an uncertainty set around the prediction. However, this approach does not apply in our setting as it requires observations of the unknown parameters, which we do not have access to. Our paper presents the first method to calibrating uncertainty sets using decision data.

\textbf{Estimate, then optimize.} Conformal IO belongs to a family of data-driven optimization methods called ``estimate, then optimize'' \citep{elmachtoub2023estimate}. Recent research suggests that even small estimation errors may be amplified in the optimization step, resulting in significant decision errors. This issue can be mitigated by training the estimation model with decision-aware losses \citep{wilder2019melding, mandi2022decision} and robustifying the optimization model \citep{chan2023got}. We take a similar approach as the second stream, yet deviate from them by i) utilizing decision data instead of observations of the unknown parameters, and ii) focusing on both the ground-truth and perceived solution quality, the latter of which has not been studied in this stream of literature.

\textbf{Preference learning.} Preference learning has been studied in the context of reinforcement learning and has recently attracted significant attention due to its application in AI alignment \citep{ji2023ai}. Existing methods focus on learning a reward function/decision policy that is maximally consistent with expert decision trajectories \citep{ng2000algorithms, wu2024inverse} or labeled preferences in the form of pair-wise comparison/ranking \citep{wirth2017survey, christiano2017deep, rafailov2024direct}. We enrich this stream of literature in two ways. First, we leverage \textit{unlabeled} decision data that are not state-action trajectories, but solutions to an optimization model. Second, instead of learning a policy that imitates expert behaviors, we aim to extract common wisdom from decision data crowd-sourced from DMs who are not necessarily experts to encourage algorithm adoption.


\section{Preliminaries} \label{sec:prelim}

In this section, we first present the problem setup (Section \ref{subsect:setup}) and then describe the challenges with the classic IO pipeline (Section \ref{subsec:risk_neutral_pipeline}). Finally, we provide intuition on why robustifying the IO pipeline would help (Section \ref{subsec:risk_averse_pipeline}). 

\subsection{Problem Setup} \label{subsect:setup}

\textbf{Data generation.} Consider a \textit{forward optimization} problem
\begin{align}
    \mathbf{FO}(\bstheta, \bfu): 
    \underset{\bfx\in \mX(\bfu)}{\textrm{minimize}}\ f(\bstheta, \bfx)
\end{align}
where $\bfx \in \bbR^n$ is the decision vector whose feasible region $\mX(\bfu)$ is non-empty and is parameterized by exogenous parameters $\bfu \in \bbR^m$, $\bstheta\in \bbR^d$ is a parameter vector, and $f: \bbR^{n \times d} \rightarrow \bbR$ is the objective function. Suppose $\bfu$ is distributed according to $\bbP_\bfu$ supported on $\mU$. There exists a ground-truth parameter vector $\bstheta^*$ that is unknown to the DM. Instead, the DM obtains a decision $\hat{\bfx}$ by solving $\mathbf{FO}(\hat{\bstheta}, \bfu)$ where $\hat{\bstheta}$ is a noisy perception of $\bstheta^*$. We assume that, while the distribution $\bbP_{\bstheta}$ of $\hat{\bstheta}$ is unknown, it is supported on a known bounded set $\bsTheta \subset \mathbb{R}^d$ and that $\bstheta^* \in \bsTheta$. Let $\bbP_{(\bstheta, \bfu)}$ denote the joint distribution of $\hat{\bstheta}$ and $\bfu$, $\tilde{\bfx}: \bsTheta \times \mU \rightarrow \bbR^n$ be an oracle that returns an optimal solution to $\mathbf{FO}$ drawn uniformly at random from $\mX^\textrm{OPT}(\bstheta, \bfu) := \argmin\left\{f(\bstheta, \bfx) \,\middle|\, \bfx \in \mX(\bfu) \right\}$. 

\textbf{Objective function.} In this paper, we focus on cases where $f$ is linear in $\bstheta$ and convex in $\bfx$. That is
$
    f(\bstheta, \bfx) = 
    \sum_{i\in [d]}
        \theta_i f_i(\bfx)
$
where $f_i: \bbR^n \rightarrow \bbR$ for all $i\in [d]$ represent some convex basis functions. This objective function generalizes the linear objective function, i.e., $f(\bstheta, \bfx) = \bstheta^\intercal \bfx$. Moreover, $\mathbf{FO}$ with this objective function can be interpreted as a multi-objective optimization model, which has been applied to model routing preferences \citep{ronnqvist2017calibrated}, radiation therapy treatment planning \citep{chan2014generalized}, and portfolio optimization \citep{dong2021wasserstein}.

\textbf{Learning task.} Given a dataset of $N$ decision-exogenous parameter pairs $\mD = \left\{ \hat{\bfx}_k, \bfu_k \right\}_{k=1}^N$, we are interested in finding a decision policy $\bar{\bfx}: \mU \rightarrow \bbR^n$ to suggest decisions for future $\bfu$. We require $\bar{\bfx}(\bfu) \in \mX(\bfu)$ for any $\bfu\in \mU$. As discussed later, $\bar{\bfx}(\bfu)$ is usually generated by solving an optimization model that may have multiple optimal solutions. So we consider randomized policies (e.g., uniformly sample from a set of optimal solutions). This is nonrestrictive because a deterministic policy can be recovered from a randomized policy that samples the deterministic solution with probability one. 

\subsubsection{Assumptions} \label{subsubsec:asm}


\begin{assumption}[I.I.D. Samples]
\label{asm:iid}
    The dataset $\mD$ is generated using $\hat{\bfx}_k := \tilde{\bfx}\left(\hat{\bstheta}_k, \bfu_k\right)$ where $\left(\hat{\bstheta}_k, \bfu_k\right)$ are i.i.d. samples from $\bbP_{(\bstheta, \bfu)}$ for all $k\in [N]$.
\end{assumption}

\begin{assumption}[Bounded Inverse Feasible Set]
\label{asm:bounded_inv_feas}
    There exists a constant $\eta\in \bbR_+$ such that, for any $\bstheta, \bstheta' \in \bsTheta^\mathrm{OPT}\left(\hat{\bfx}, \bfu\right)$, for some $\hat{\bfx} \in \mX(\bfu)$ and $\bfu \in \mU$, we have $\|\bstheta - \bstheta'\|_2 \leq \eta$, where
    \begin{equation}
        \bsTheta^\mathrm{OPT}(\bfx, \bfu) 
        := 
        \left\{ 
            \bstheta\in \bbR^d 
        \,\middle|\, 
            \bfx \in \mX^\mathrm{OPT}(\bstheta, \bfu), \|\bstheta\|_2 = 1 
        \right\}.  
    \end{equation}
\end{assumption}

\begin{assumption}[Bounded Divergence] \label{asm:bnd_div}
    There exists a constant $\sigma \in \bbR_+$ such that $\|\bbE(\hat{\bstheta}) - \bstheta^*\|_2 \leq \sigma$. 
\end{assumption}

Assumption \ref{asm:iid} is standard in the ML and IO literature. Assumption \ref{asm:bounded_inv_feas} is mild because $\bsTheta^\textrm{OPT}\left(\hat{\bfx}, \bfu\right)$ is by definition bounded and is usually much smaller than $\bsTheta$. Assumption \ref{asm:bnd_div} states that the $l_2$ distance between the expected perceived parameters and the ground-truth parameters is upper bounded. It is reasonable in many real-world settings. For example, a driver's perceived travel cost ($\hat{\bstheta}$) should not be too different from the travel time ($\bstheta^*$) as the latter is an important factor that drivers consider.

\subsubsection{Evaluation Metrics}

\begin{definition}
   The actual optimality gap (AOG) of a decision policy $\bar{\bfx}$ is defined as 
    \begin{equation} \label{eq:aog}
        \mathrm{AOG}(\bar{\bfx})
        := 
        \bbE
        \left[
            f\left(\bstheta^*, \bar{\bfx}(\bfu) \right) -  f\left(\bstheta^*, \tilde{\bfx}(\bstheta^*, \bfu) \right)
        \right]
    \end{equation}
    where the expectation is taken over the joint distribution of the random variable $\bfu$ and the decision sampled using the possibly randomized policy $\bar{\bfx}$.
\end{definition}

\begin{definition}
    The perceived optimality gap (POG) of a decision policy $\bar{\bfx}$ is defined as 
    \begin{equation} \label{eq:pog}
        \mathrm{POG}(\bar{\bfx})
        := 
        \bbE
        \left[
            f\left(\hat{\bstheta}, \bar{\bfx}(\bfu) \right) 
            - f\left(\hat{\bstheta}, \tilde{\bfx}\left(\hat{\bstheta}, \bfu\right) \right)
        \right].
    \end{equation}
    where the expectation is taken with respect to the randomness in  $\hat{\bstheta}$, $\bfu$, and possibly $\bar{\bfx}$. 
\end{definition}

AOG is an objective performance measure. Achieving a low AOG means that $\bar{\bfx}$ can generate high-quality decisions. In contrast, POG is a subjective measure that depends on the DM's perception of the problem. Achieving a low POG is critical to mitigate algorithm aversion \citep{burton2020systematic}. 

\subsection{An Inverse Optimization Pipeline} \label{subsec:risk_neutral_pipeline}
Finding $\bar{\bfx}$ is challenging for three reasons. First, unlike many ML tasks where the prediction target is unconstrained, we require $\bar{\bfx}(\bfu)$ to be feasible to $\mathbf{FO}$ which may involve a large number of constraints. Supervised learning approaches that predict $\hat{\bfx}$ based on $\bfu$ can often fail as they typically do not provide feasibility guarantees. An optimization module is often needed to recover feasibility or produce feasible solutions based on $\bfu$ and some estimated $\bar{\bstheta}$. Second, we do not directly observe $\bstheta^*$ or $\hat{\bstheta}$, which precludes using classic ML techniques to estimate them. Finally, AOG and POG may not necessarily align with each other, so we are essentially dealing with a bi-objective problem.

\begin{figure}[th]
    \centering
    \includegraphics[width=.85\textwidth]{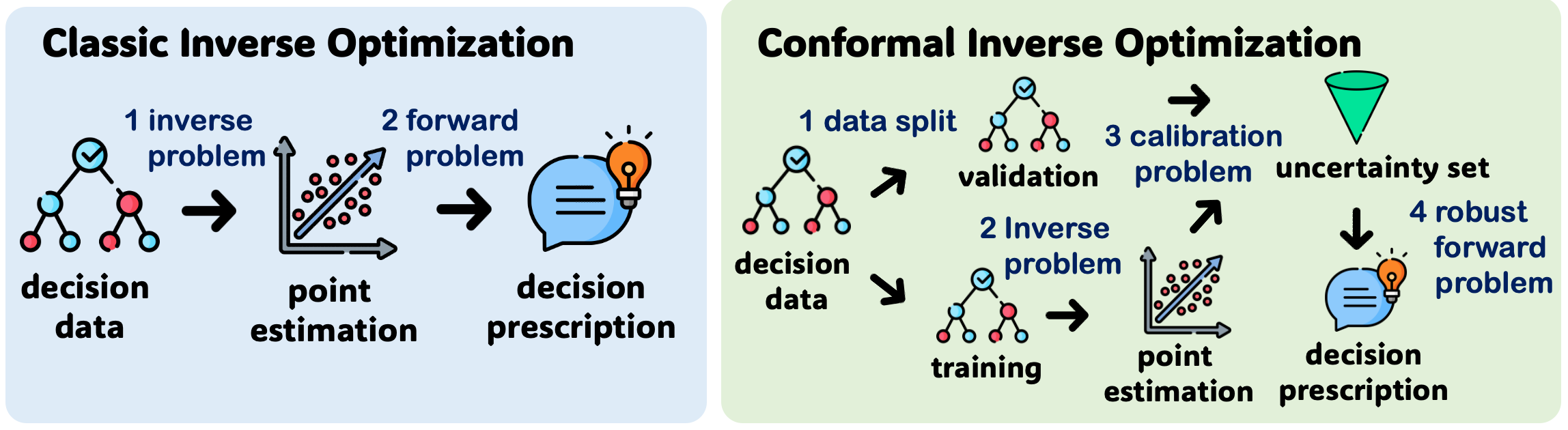}
    \caption{Classic and conformal IO pipelines.}
    \label{fig:pipelines}
\end{figure}

In light of the first two challenges, a classic IO pipeline (visualized in Figure \ref{fig:pipelines}) has been proposed to first obtain a point estimation $\bar{\bstheta}$ of the unknown parameters and then employ a policy $\bar{\bfx}_\textrm{IO}(\bfu) := \tilde{\bfx}(\bar{\bstheta}, \bfu)$ to prescribe decisions for any $\bfu\in \mU$ \citep{ronnqvist2017calibrated, babier2020knowledge}.
Specifically, we can estimate the parameters by solving the following \textit{inverse optimization} problem
\begin{align}
     \mathbf{IO}(\mD): \underset{\bstheta \in \bsTheta}{\textrm{minimize}} \ 
    & \frac{1}{N} \sum_{k\in [N]} 
        \ell\left(\hat{\bfx}_k, \mX^\textrm{OPT}(\bstheta, \bfu_k)\right),
\end{align}
where $\ell$ is a non-negative loss function that returns 0 when $\hat{\bfx}_k \in \mX^\textrm{OPT}(\bstheta, \bfu_k)$. For instance, the following loss function is a popular choice in the literature.

\begin{definition}
    The sub-optimality loss of $\bstheta$ is given by
    \begin{equation} \label{eq:sub_opt_loss}
        \ell_\textrm{S}\left(\hat{\bfx}, \mX^\textrm{OPT}(\bstheta, \bfu)\right)
        := \max_{\bfx\in \mX(\bfu)} f(\bstheta, \hat{\bfx}) - f(\bstheta, \bfx).
    \end{equation}
\end{definition}

\begin{figure}
    \centering
    \includegraphics[width=0.50\textwidth]{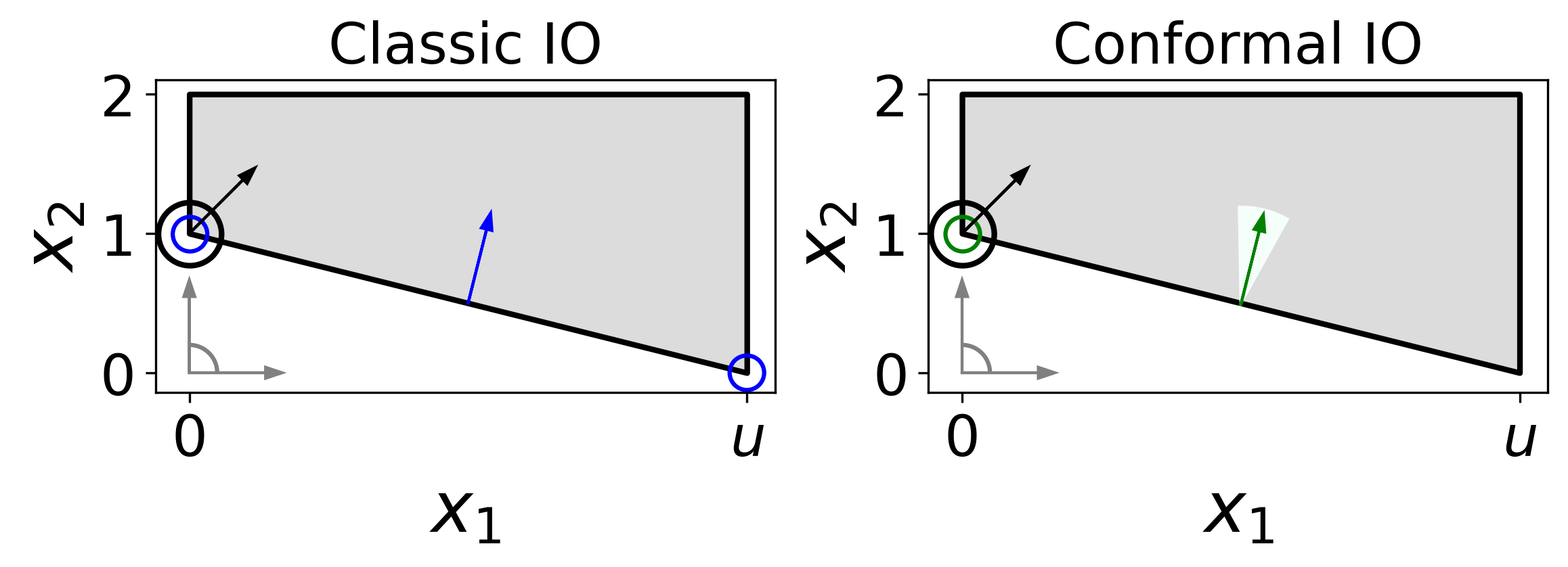}
    \caption{Illustration of Example \ref{ex:2d}. The gray areas are the feasible region $\mX(u)$. The black arrows are the ground-truth parameter $\bstheta^*$. The gray arrows are the extreme rays of $\bsTheta$. The blue and green arrows are the point estimation $\bar{\bstheta}$. The green area is the uncertainty set $\mC(\bar{\bstheta}, \alpha)$. The black circles are the optimal solution to $\mathbf{FO}(\bstheta^*, u)$. The blue and green circles are the suggested decisions. Note that $\bar{\bfx}_\textrm{IO}$ may suggest any decisions on the facet of $x_1 + ux_2 \geq u$, which are omitted for clarity.
    }
    \label{fig:2d_exp}
\end{figure}

The sub-optimality loss penalizes the optimality gap achieved by the observed decision under the estimated parameters. As remarked by \citet{mohajerin2018dataIO}, this loss function has better computational properties than its alternatives as it is convex in the unknown parameters. In fact, when the dataset $\mD$ is large in size and the unknown parameters $\bstheta$ are high-dimensional, the sub-optimality loss is usually the only loss function that leads to a tractable inverse problem, although it does not enjoy properties such as statistical consistency (see \citet{chan2023inverse} for detailed discussions). While such a trade-off is acceptable in some applications, we suggest that it is undesirable in our setting because the resulting policy can achieve arbitrarily large AOG and POG. To see this, consider the following example that satisfies Assumptions \ref{asm:iid}--\ref{asm:bnd_div}. This example is visualized in Figure \ref{fig:2d_exp}. 

\begin{example} \label{ex:2d}
    Let $\mathbf{FO}(\theta, u)$ be the following problem
    \begin{subequations}
    \label{prob:2d_exp}
    \begin{align}
        \mathrm{minimize} \quad 
            & \theta_1 x_1  + \theta_2 x_2 \\
        \mathrm{subject\ \ to} \quad
            & x_1 + u x_2 \geq u \\
            & 0 \leq x_1 \leq u \\
            & 0 \leq x_2 \leq 2.
    \end{align}
    \end{subequations}
    Let the ground-truth $\bstheta^* = \left(\cos(\pi/4), \sin(\pi/4) \right)$ and $\mU = \{u\}$ where $u > 1$ is a real constant. We are given a dataset $\mD = \left\{\hat{\bfx}_k, u \right\}_{k=1}^N$ where $\hat{\bfx}_k = \tilde{\bfx}(\hat{\bstheta}_k, u)$ with $\hat{\bstheta}_k$ uniformly and independently drawn from $\bsTheta = \left\{(\cos\delta,\, \sin\delta) \,\middle|\, \delta \in (0, \pi/2) \right\}$ for all $k\in [N]$.
\end{example}

\begin{lemma} \label{lemma:2d_exp_asy_est}
    In Example \ref{ex:2d}, let $\bar{\bstheta}_N$ denote an optimal solution to $\mathbf{IO}(\mD)$ with the sub-optimality loss~\eqref{eq:sub_opt_loss}, we have $\bbP\left(
            \bar{\bstheta}_N = \bstheta_u
        \right) 
        \rightarrow 1$ 
    as $N \rightarrow \infty$, where $\bstheta_u := \left(1/ \sqrt{1 + u^2}, u/ \sqrt{1 + u^2} \right)$.
\end{lemma}

Lemma \ref{lemma:2d_exp_asy_est} shows that, when using $\mathbf{IO}(\mD)$ with the sub-optimality loss to estimate the unknown parameter in Example \ref{ex:2d}, the probability of the estimated parameter being $\bstheta_u$ converges to one asymptotically. This implies that asymptotically we are almost certain that $\bar{\bfx}_\mathrm{IO}(u) = \tilde{\bfx}(\bstheta_u, u)$, i.e. the policy that samples uniformly from the facet corresponding to the constraint $x_1 + ux_2 \geq u$.
As a result, $\bar{\bfx}_\textrm{IO}$ can achieve arbitrarily large AOG and POG when $u$ is set to a large enough value. 

\begin{proposition} \label{prop:unbd_pog_aog}
    In Example \ref{ex:2d}, let $\bar{\bfx}_\mathrm{IO}(u) = \tilde{\bfx}(\bstheta_u, u)$. For any $v\in \bbR_+$ there exists some $\bar{u} > 1$ such that $\mathrm{AOG}(\bar{\bfx}_\mathrm{IO}) > v$ and $\mathrm{POG}(\bar{\bfx}_\mathrm{IO}) > v$ for any $u > \bar{u}$.
\end{proposition}


\subsection{Robustifying the Inverse Optimization Pipeline} \label{subsec:risk_averse_pipeline}

A natural idea to improve the AOG and POG of $\bar{\bfx}$ is to robustify the decision pipeline. Specifically, instead of solving $\mathbf{FO}$ with some estimated parameters $\bar{\bstheta}$, we solve the following \textit{robust forward optimization problem} with an uncertainty set around $\bar{\bstheta}$ to prescribe decisions (final steps in Figure \ref{fig:pipelines}). 
\begin{align}
    \mathbf{RFO}\left(\mC(\bar{\bstheta}, \bsalpha), \bfu\right): 
    \underset{\bfx\in \mX(\bfu)}{\textrm{minimize}} \ 
    \underset{\bstheta\in \mC(\bar{\bstheta}, \bsalpha)}{\textrm{maximize}} 
    \ f(\bstheta, \bfx)
\end{align}
where $\mC$ is an uncertainty set with $\bar{\bstheta}$ being its center and $\bsalpha$ representing parameters that control its shape/size. In this paper, since the objective function $f$ is linear in $\bstheta$, the optimal solution to $\mathbf{FO}$ is invariant to the scale of $\bstheta$, i.e., if $\bfx \in \mX^\textrm{OPT}(\bstheta, \bfu)$, then $\bfx\in \mX^\textrm{OPT}(\beta\bstheta, \bfu)$ for any $\beta \in \bbR_+$. We thus focus on the following uncertainty set without loss of generality.
\begin{equation} \label{def:cone}
    \mC(\bar{\bstheta}, \alpha) := 
    \left\{
        \bstheta\in \bbR^d
    \,|\, 
        \|\bstheta\|_2 = 1,\, \bstheta^\intercal \bar{\bstheta} \geq \cos\alpha
    \right\}
\end{equation}
where $\alpha \in (0, \pi]$ represents the max angle between $\bar{\bstheta}$ and any vector in the uncertainty set. 

\begin{remark}
    An alternative approach to robustify the classic IO pipeline is to replace $\mathbf{IO}$ with a distributionally robust IO for parameter estimation \citep{mohajerin2018dataIO}. However, this approach would not help as $\bstheta_u$ is still an optimal solution to the distributionally robust IO in Example \ref{ex:2d}, so the resulting AOG and POG can still be arbitrarily large. See Appendix \ref{sec:alter_robust} for the complete statement and discussions.
\end{remark}

Now, in Example \ref{ex:2d}, we analyze the performance of a policy that utilizes $\mathbf{RFO}$ to prescribe decisions.




\begin{lemma}
\label{lemma:2d_robust_aog_pog}
    In Example \ref{ex:2d}, let $\bar{\bfx}_\mathrm{CIO}(u)$ be an optimal solution to $\mathbf{RFO} \left(\mC(\bar{\bstheta}_N, \alpha), u \right)$ where $\bar{\bstheta}_N$ is an optimal solution to $\mathbf{IO}(\mD)$ with the sub-optimality loss~\eqref{eq:sub_opt_loss}. When $\alpha \in (0, \pi/2)$, we have $\bbP\left[\mathrm{AOG}(\bfx_\mathrm{CIO}) = 0\right] \rightarrow 1$ and $\bbP\left[\mathrm{POG}(\bar{\bfx}_\mathrm{CIO}) < \pi/2\sqrt{2}\right] \rightarrow 1$ as $N\rightarrow \infty$.
\end{lemma}


Lemmas \ref{lemma:2d_robust_aog_pog} shows that, when using $\mathbf{RFO}$ to prescribe new decisions, the probability of achieving upper-bounded AOG and POG converges to one as $N$ goes to infinity, as long as $\alpha\in (0,\pi/2)$. These bounds are independent of $u$, in contrast to the AOG and POG of classic IO that can be arbitrarily large as $u$ changes (Proposition \ref{prop:unbd_pog_aog}). However, the performance of this approach still depends on the choice of $\alpha$, which is non-trivial when $\mathbf{FO}$ is more complex than a two-dimensional linear program. 


\section{Conformal Inverse Optimization} \label{sec:cio}

In this section, we present a principled approach to learn uncertainty sets that lead to provable performance guarantees. As presented later, the learned uncertainty set contains parameters that make the next DM's decision optimal with a specified probability. We call this approach conformal IO due to its connection to conformal prediction \citep{vovk2005algorithmic}, which aims to predict a set that contains the next prediction target with a specified probability. As illustrated in Figure \ref{fig:pipelines}, conformal IO has three steps: i) data split, ii) point estimation, and iii) uncertainty set calibration. We present these steps in Section \ref{subsec:cio_method} and analyze the properties of conformal IO in Section \ref{subsec:prop}.

\subsection{Learning an Uncertainty Set} \label{subsec:cio_method}

\textbf{Data split.} We first split the dataset $\mD$ into training and validation sets, namely $\mD_\textrm{train}$ and $\mD_\textrm{val}$. Let $\mK_\textrm{train}$ and $\mK_\textrm{val}$ index $\mD_\textrm{train}$ and $\mD_\textrm{val}$, respectively, while $N_\textrm{train} = |\mD_\textrm{train}|$ and $N_\textrm{val} = |\mD_\textrm{val}|$. 

\textbf{Point estimation.} Given a training set $\mD_\textrm{train}$, we apply data-driven IO techniques to obtain a point estimation $\bar{\bstheta}$ of the unknown parameters. The most straightforward way is to solve $\mathbf{IO}(\mD_\textrm{train})$ with any loss function. Alternatively, one may consider using end-to-end learning and optimization methods that do not require observations of the parameter vectors, e.g., the one proposed by \citet{berthet2020learn}. The point estimation can also come from other sources, e.g., from an ML model that predicts the parameters. Our calibration method is independent of the point estimation method.

\textbf{Uncertainty set calibration.} Given a point estimation $\bar{\bstheta}$, we calibrate an uncertainty set that, with a specified probability, contains parameters that make the next observed decision optimal. This property is critical for the results in Section \ref{subsec:prop} to hold. While we can naively achieve this by setting $\alpha = \pi$, the resulting $\mathbf{RFO}$ may generate overly conservative decisions. Hence, we are interested in learning the smallest uncertainty set that satisfies this condition. We solve the following \textit{calibration problem} 
\begin{subequations}
\label{prob:calib}
\begin{align}
    \label{cp_obj:min_angle}
    \mathbf{CP}(\bar{\bstheta}, \mD_\textrm{val}, \gamma): 
    \underset{\alpha, \{\bstheta_k\}_{k\in \mK_\textrm{val}}}
    {\textrm{minimize}} \quad
    & \alpha 
    \\
    \textrm{subject to} \quad 
    \label{cp_con:inv_feas}
    & \hat{\bfx}_k \in \mX^\mathrm{OPT}(\bstheta_k, \bfu_k),
        \  \forall k\in \mK_\mathrm{val} \\
    \label{cp_con:probability}
    & \sum_{k\in \mK_\textrm{val}}
        \mathds{1}
        \left[
            \bstheta_k \in  
            \mC(\bar{\bstheta}, \alpha) 
        \right]
    \geq \gamma (N_\textrm{val} + 1) \\
    \label{cp_con:reg}
    & \|\bstheta_k \|_2 = 1,
        \  \forall k\in \mK_\mathrm{val}\\
    \label{cp_con:alpha_range}
    & 0 \leq \alpha \leq \pi,
\end{align}
\end{subequations}
where decision $\alpha$ controls the size of the uncertainty set, $\bstheta_k$ represent a possible parameter vector associated with data point $k\in \mK_\textrm{val}$, $\gamma \in [0, 1]$ is a DM-specified confidence level. Constraints~\eqref{cp_con:inv_feas} ensure that $\bstheta_k$ can make the decision $\hat{\bfx}_k$ optimal for $k\in \mK_\textrm{val}$. Constraint~\eqref{cp_con:probability} ensures that at least $\gamma(N_\textrm{val} + 1)$ of the decisions in $\mD_\textrm{val}$ can find a vector in $\mC$ that makes it optimal. Constraints \eqref{cp_con:reg} ensure that the parameter vectors are on the unit sphere. 

\begin{remark}[Optimality Conditions]
    The specific form of Constraints~\eqref{cp_con:inv_feas} depends on the structure of $\mathbf{FO}$. For example, when $\mathbf{FO}$ is a linear program, Constraints~\eqref{cp_con:inv_feas} can be replaced with the dual feasibility and strong duality constraints. When the $\mathbf{FO}$ is a general convex optimization problem, we can use the KKT conditions. For non-convex forward problems, we can replace Constraints~\eqref{cp_con:inv_feas} with 
    $f(\bstheta_k, \hat{\bfx}_k) \leq f(\bstheta_k, \bfx)$ for all $\bfx \in \mX(\bfu)$, which can be generated in a cutting-plane fashion.
\end{remark}

\begin{remark}[Feasibility]
    For $\mathbf{CP}$ to be feasible, we require, for each observed decision, there exists a $\bstheta\in \bsTheta$ that make it optimal. This condition holds for a range of problems, e.g., routing problems and the knapsack problem, even if the DM is subject to bounded rationality, i.e., the DM settles for suboptimal solutions due to cognitive/computational limitations. For problems where this condition is violated, we may pre-process $\mD_\textrm{val}$ to project $\hat{\bfx}$ to a point in $\mX(\bfu)$ such that the condition is satisfied.
\end{remark}

Solving $\mathbf{CP}$ is hard. First, $\mathbf{CP}$ is non-convex due to Constraints~\eqref{cp_con:reg}. Second, Constraints~\eqref{cp_con:inv_feas} involve the optimality conditions of $N_\textrm{val}$ problems, so the size of $\mathbf{CP}$ grows quickly as $N_\textrm{val}$ increases. Nevertheless, considering a large $\mD_\textrm{val}$ is critical to ensure desirable properties of the learned uncertainty set (Section \ref{subsec:prop}). Below we introduce a decomposition method to solve $\mathbf{CP}$ efficiently. 

\begin{theorem}
    \label{thm:cali_reform_cone}
    Let $\mD_\mathrm{val}$ be a dataset, $\gamma\in [0, 1]$,  $\bar{\bstheta}\in \bbR^d$, $\tau = \left\lceil \gamma (N_\mathrm{val} + 1) \right\rceil$ and $\Gamma_\tau$ be an operator that returns the $\tau^\mathrm{th}$ largest value in a set. The optimal solution to $\mathbf{CP}(\bar{\bstheta}, \mD_\mathrm{val}, \gamma)$ is 
    $\alpha_\gamma := \arccos\left(\Gamma_\tau\left(\{c_k\}_{k\in \mK_\textrm{val}}\right) \right)$ 
    with $c_k := 
    \max_{\bstheta_k}
    \left\{
        \bstheta_k^\intercal \bar{\bstheta} 
    \,\middle|\,
        \hat{\bfx}_k \in \mX^\mathrm{OPT}(\bstheta_k, \bfu_k), \, \|\bstheta_k\|_2 \leq 1
    \right\}
    $.
\end{theorem}

Theorem \ref{thm:cali_reform_cone} states that we can solve $\mathbf{CP}$ by first solving $N_\textrm{val}$ optimization problems whose size is independent of $N_\textrm{val}$ and then find a quantile in a set of $N_\textrm{val}$ elements. The first step is parallelizable and the second step can be done in $O\left(N_\textrm{val}\log(\tau)\right)$ time. Since the problem required for evaluating $c_k$ is a maximization problem, we can replace the constraint $\|\bstheta_k\|_2 = 1$ with $\|\bstheta_k\|_2 \leq 1$, so this problem is convex when the forward problem is a linear program. 



\subsection{Properties of Conformal IO} \label{subsec:prop}

\begin{theorem}[Uncertainty Set Validity]
    \label{thm:valid}
    Let $\mD_\mathrm{val}$ be a dataset that satisfies Assumption \ref{asm:iid}, $(\hat{\bstheta}, \bfu)$ be a new i.i.d. sample from $\bbP_{(\bstheta, \bfu)}$, $\hat{\bfx} = \tilde{\bfx}(\hat{\bstheta}, \bfu)$, $\hat{\bsTheta} := \bsTheta^\mathrm{OPT}(\hat{\bfx}, \bfu)$, and $\alpha_\gamma$ be an optimal solution to $\mathbf{CP}(\bar{\bstheta}, \mD_\mathrm{val}, \gamma)$ where $\bar{\bstheta}\in \bbR^d$. We have, for any $\gamma\in [0,  N_\textrm{val}/(N_\textrm{val} + 1)]$, 
    \begin{equation} \label{ineq:conserv_valid}
    \bbP
    \left(
        \hat{\bsTheta} \cap \mC(\bar{\bstheta}, \alpha_\gamma)
        \neq
        \emptyset
    \right)
    \geq \gamma.
\end{equation}
For any $\gamma\in [0, 1]$, with probability at least $1 - 1/N_\mathrm{val}$,
\begin{equation} \label{ineq:asy_exact_valid}
    \left|
    \bbP
    \left(
        \hat{\bsTheta} \cap \mC(\bar{\bstheta}, \alpha_\gamma)
        \neq
        \emptyset
    \right)
    - \gamma
    \right|
    \leq
    \sqrt{\frac{8\log(N_\textrm{val} + 1) + 2\log N_\textrm{val}}{N_\textrm{val}}}
    + \frac{2}{N_\textrm{val}}.
\end{equation} 
\end{theorem}

Theorem \ref{thm:valid} states that our learned uncertainty set is conservatively valid and asymptotically exact \citep{vovk2005algorithmic}. More specifically, first, our method will produce a set that contains a $\bstheta$ that makes the next DM's decision optimal no less than $\gamma$ of the time that it is used (conservatively valid). The probability in Inequality~\eqref{ineq:conserv_valid} is with respect to the joint distribution over $\mD_\textrm{val}$ and the new sample. Second, once the set is given, we have high confidence that, the probability of the next DM's decision being covered is within $\epsilon(N_\textrm{val})$ from $\gamma$. The probability in Inequality~\eqref{ineq:asy_exact_valid} is with respect to the new sample, while the high confidence is with respect to the draw of the validation data set. Overall, we have the almost sure convergence of $\bbP
    (
        \hat{\bsTheta} \cap \mC(\bar{\bstheta}, \alpha_\gamma)
        \neq
        \emptyset
    )$ 
to $\gamma$ as $N$ goes to infinity.

Now, we relate the validity results to the performance of conformal IO. The following Lemma is an immediate result of the objective function $f$ being linear in $\bstheta$.

\begin{lemma}\label{lemma:lipschiz}
    For any $\hat{\bfx} \in \tilde{\bfx}\left(\hat{\bstheta}, \bfu \right)$ and $\left(\hat{\bstheta}, \bfu\right) \in \bsTheta \times \mU$, there exists a constant $\nu\left(\hat{\bfx} \right) \in \bbR_+$ such that, for any $\bstheta, \bstheta' \in \bsTheta$, we have $f\left(\bstheta, \hat{\bfx}\right) - f\left(\bstheta', \hat{\bfx}\right) \leq \nu\left(\hat{\bfx} \right) \|\bstheta - \bstheta'\|_2$.
\end{lemma}

\begin{theorem}[POG Bound] \label{thm:pog}
    Let $\bar{\bfx}_\mathrm{CIO}(\bfu)$ be an optimal solution to $\mathbf{RFO}\left(\mC(\bar{\bstheta}, \alpha_1), \bfu \right)$ for any $\bfu\in \mU$, where $\bar{\bstheta} \in \bbR^d$ and $\alpha_1$ are chosen such that, for a new sample $(\bstheta', \bfu')$ from $\bbP_{(\bstheta, \bfu)}$ and $\bfx' = \tilde{\bfx}(\bstheta', \bfu')$, $\bbP \left( \mC(\bar{\bstheta}, \alpha_1) \cap \bsTheta^\mathrm{OPT}(\bfu', \bfx') \neq \emptyset \right) = 1$. If Assumptions \ref{asm:bounded_inv_feas}--\ref{asm:bnd_div} hold, then  
    \begin{equation}
        \mathrm{POG}(\bar{\bfx}_\mathrm{CIO}) 
        \leq 
        (\eta - 2\cos2\alpha_1 + 2)\mu
        + \eta \mu_\mathrm{CIO},
    \end{equation}
    and 
    \begin{equation}
        \mathrm{AOG}(\bar{\bfx}_\mathrm{CIO}) 
        \leq 
        (2 - 2 \cos 2\alpha_1 + \eta + \sigma) \mu^*   
        + (\eta + \sigma) \mu_\mathrm{CIO},
    \end{equation}
    where $\mu := \bbE [\nu(\tilde{\bfx}(\hat{\bstheta}, \bfu))]$, $\mu_\mathrm{CIO} := \bbE (\nu[ \bar{\bfx}_\mathrm{CIO}(\bfu) ])$, and $\mu^* := \bbE\left(\nu[ \tilde{\bfx}(\bstheta^*, \bfu)] \right)$.
\end{theorem}


Theorem \ref{thm:pog} states that, when the uncertainty set contains a $\bstheta$ that makes the next DM's decision optimal almost surely, conformal IO achieves upper-bounded POG and AOG. While one can meet this condition by using a large $\alpha$, the resulting bounds would be inefficient as they increase as $\alpha$ increases, reflecting that the prescribed decisions are overly conservative. Instead, we can use $\mathbf{CP}$ to obtain an $\alpha$ that achieves close-to-100\% coverage. Our bounds involve problem-specific constants. To demonstrate tightness, we present our bounds in Example \ref{ex:2d} with $\alpha=\pi/4$ in Table \ref{tab:bound_vals}. Our bounds closely follow the performance of conformal IO, which outperforms classic IO by a large margin.

\begin{table}[!ht]
\centering \footnotesize
\caption{Performance profile of classic and conformal IO in Example \ref{ex:2d}.} \label{tab:bound_vals}
\begin{tabular}{@{}lrrrrrrrrr@{}}
\toprule
                   & \multicolumn{4}{c}{AOG}  &   & \multicolumn{4}{c}{POG}     \\ 
\cmidrule(lr){2-5} \cmidrule(l){7-10}
$u$                & 2    & 10   & 50    & 100   &   & 2    & 10   & 50    & 100   \\ \midrule
Classic IO         & 0.35 & 3.18 & 17.32 & 35    &   & 0.74 & 4.55 & 24.51 & 49.50 \\
Conformal IO       & 0.00 & 0.00 & 0.00  & 0.00  &   & 0.70 & 0.16 & 0.03  & 0.02  \\
Conformal IO bound & 0.70 & 0.05 & 0.002 & 0.001 &   & 1.58 & 0.67 & 0.15  & 0.08  \\ \bottomrule
\end{tabular}
\end{table}

\section{Numerical Studies} \label{sec:comp}


\textbf{Data generation.} We consider two forward problems: The shortest path problem (linear program, $d = 120$) and knapsack problem (integer program, $d = 10$). See Appendix \ref{subsec:fwd_prob} for their formulations. We use synthetic instances, which is a common practice as there is no well established IO benchmark \citep{tan2020learn, dong2018generalized}. For both problems, we randomly generate a ground-truth parameters $\bstheta^*$ and a dataset of $N=1000$ DMs. For each DM $k\in [N]$, we generate her perceived parameters as $\hat{\theta}^i_k = (\theta^{i*}_k * p^i_k + \epsilon^i_k)^+ + \epsilon_0$ for $i\in [d]$ where $p^i_k$ is uniformly drawn from $[1/2, 2]$, $\epsilon^i_k$ is drawn from a normal distribution with mean 0 and standard deviation 1, and $\epsilon_0 = 0.1$. For the shortest path problem, parameters $\bfu_k$ represent a random origin-destination pair on the network. For the knapsack problem, $\bfu_k$ correspond to the weights of different items and budget of the DM. The item weight $w^i$ for $i\in [d]$ are uniformly drawn from $[1, 10]$ and are shared among DMs. For each DM $k\in [N]$, we generate a budget $u_k = q_k\sum_{i}w^i$ where $q_k$ is uniformly drawn from $[1/5, 5]$. 

\textbf{Experiment design.} Conformal IO is compatible with any approach that can provide a point estimation of unknown parameters using \textit{decision data} (Step 2 in Section \ref{subsec:cio_method}). To the best of our knowledge, i) $\mathbf{IO}$ with the sub-optimality loss and ii) the PFYL approach from \citet{berthet2020learn} are the only two methods that can perform this task \textit{at scale}. We thus implement conformal IO with these two methods. They also serve as our baselines. We call both i) and ii) classic IO to emphasize that they rely on a point estimation for decision prescription, although PFYL is not an IO approach. See Appendix \ref{ec:implementation} for implementation details. In all experiments, conformal IO uses the training set for point estimation and the validation set for calibration, while classic IO uses the union of the training and validation sets for point estimation. So, they have access to the same amount of data and are evaluated on the same test set. Unless otherwise noted, experiments are based on a 60/20/20 train-validation-test split and are repeated 10 times with different random seeds. 

\begin{figure}[!ht]
    \centering
    \includegraphics[width=.7\textwidth]{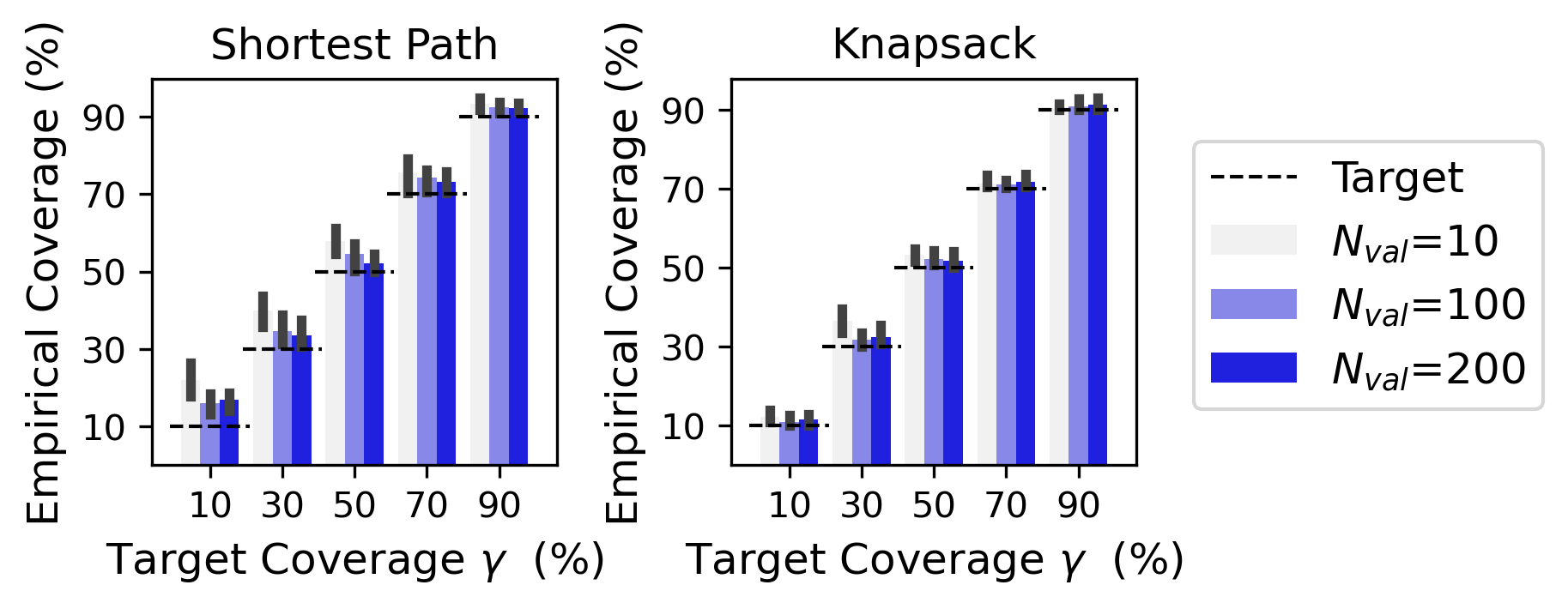}
    \caption{Empirical coverage achieved by the learned uncertainty set (error bar = range).} \label{fig:cov}
    \label{fig:enter-label}
\end{figure}

\textbf{Uncertainty validity.} To verify Theorem \ref{thm:valid}, we empirically evaluate the out-of-sample coverage achieved by our uncertainty set under different target levels $\gamma$ and sample sizes $N_\textrm{val}$. The point estimation is generated by $\mathbf{IO}$ with the sub-optimality loss. As shown in Figure \ref{fig:cov}, when the validation set is small ($N_\textrm{val} = 50$), we always achieve the specified target but $\mC(\bar{\bstheta}, \alpha_\gamma)$ tends to over-cover (conservatively valid). When using larger validation sets ($N_\textrm{val} \in \{100, 200\}$), our coverage level gets closer to the specified $\gamma$ (asymptotic exact). These empirical findings echo our theoretical analysis.

\begin{figure}[ht]
    \centering
    \includegraphics[width=0.8\textwidth]{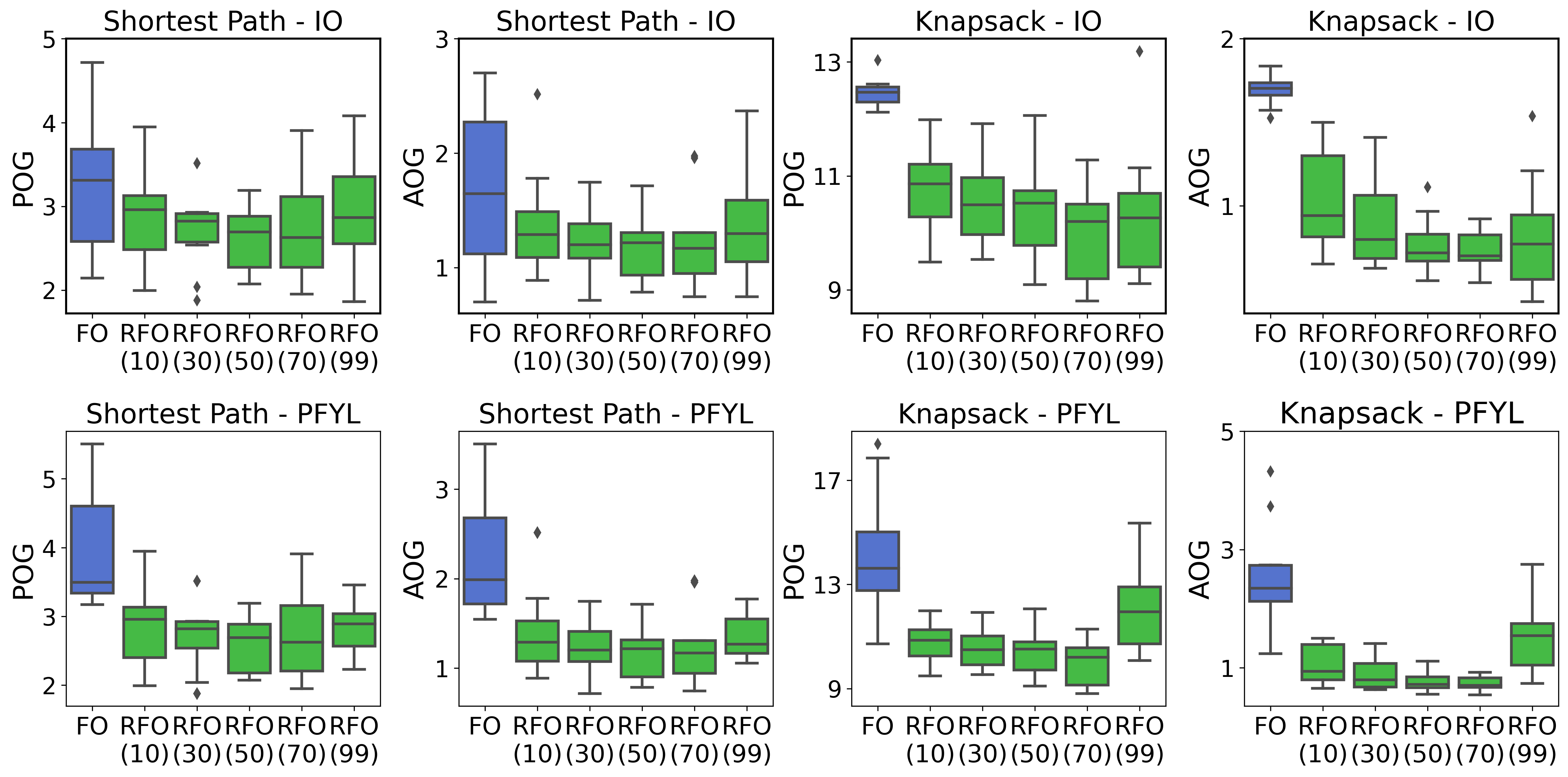}
    \caption{Performance profile of classic (blue) and conformal IO (green).}
    \label{fig:pog_aog}
\end{figure}

\textbf{The value of robustness}. As shown in Figure \ref{fig:pog_aog}, solving $\mathbf{RFO}$ with an uncertainty set learned by conformal IO leads to decisions of lower AOG and POG, compared to solving $\mathbf{FO}$ with a point estimation from classic IO. On average, when varying $\gamma$, our approach improves AOG by 20.1--30.4\% and POG by 15.0--23.2\% for the shortest path problem, and improves AOG by 40.3--57.0\% and POG by 13.5--20.1\% for the knapsack problem. The improvement is orthogonal to the point estimation method. Our decisions are of higher quality and better align with human intuition.

\textbf{Computational efficiency.} As shown in Table \ref{tab:comp_eff}, conformal IO and classic IO require similar training times. When $\mathbf{FO}$ is an integer program (knapsack), the training of conformal IO is even faster because it replaces a relatively large inverse integer program (associated with $\mD_\textrm{train} \cup \mD_\textrm{val}$), which is notoriously difficult to solve \citep{bodur2022inverse}, with a smaller inverse integer program (associated with $\mD_\textrm{train}$) and a set of small calibration problems that are parallelizable (Theorem \ref{thm:cali_reform_cone}). At the prediction time, our method achieves lower AOG and POG at the cost of solving a more challenging $\mathbf{RFO}$. Nevertheless, the solution time of $\mathbf{RFO}$ is within one second in our instances.   

\begin{table}[!ht]
\centering
\footnotesize
\caption{Average (std) computational time of classic and conformal IO in seconds.} \label{tab:comp_eff}
\begin{tabular}{@{}lrrrr@{}}
\toprule
\multicolumn{1}{c}{}        & \multicolumn{2}{c}{Training}                                      & \multicolumn{2}{c}{Prediction (per decision)}                     \\ \cmidrule(l){2-5} 
\multicolumn{1}{c}{Problem} & \multicolumn{1}{c}{Classic IO} & \multicolumn{1}{c}{Conformal IO} & \multicolumn{1}{c}{$\mathbf{FO}$} & \multicolumn{1}{c}{$\mathbf{RFO}$} \\ \midrule
Shortest Path               & 0.18 (0.02)                    & 0.27 (0.03)                      & 0.01 (0.00)                    & 0.63 (0.12)                      \\
Knapsack                    & 2.47 (0.37)                    & 1.95 (0.32)                      & 0.01 (0.00)                    & 0.44 (0.15)                      \\ \bottomrule
\end{tabular}
\end{table}

\begin{figure}[!ht]
    \centering
    \includegraphics[width=.5\textwidth]{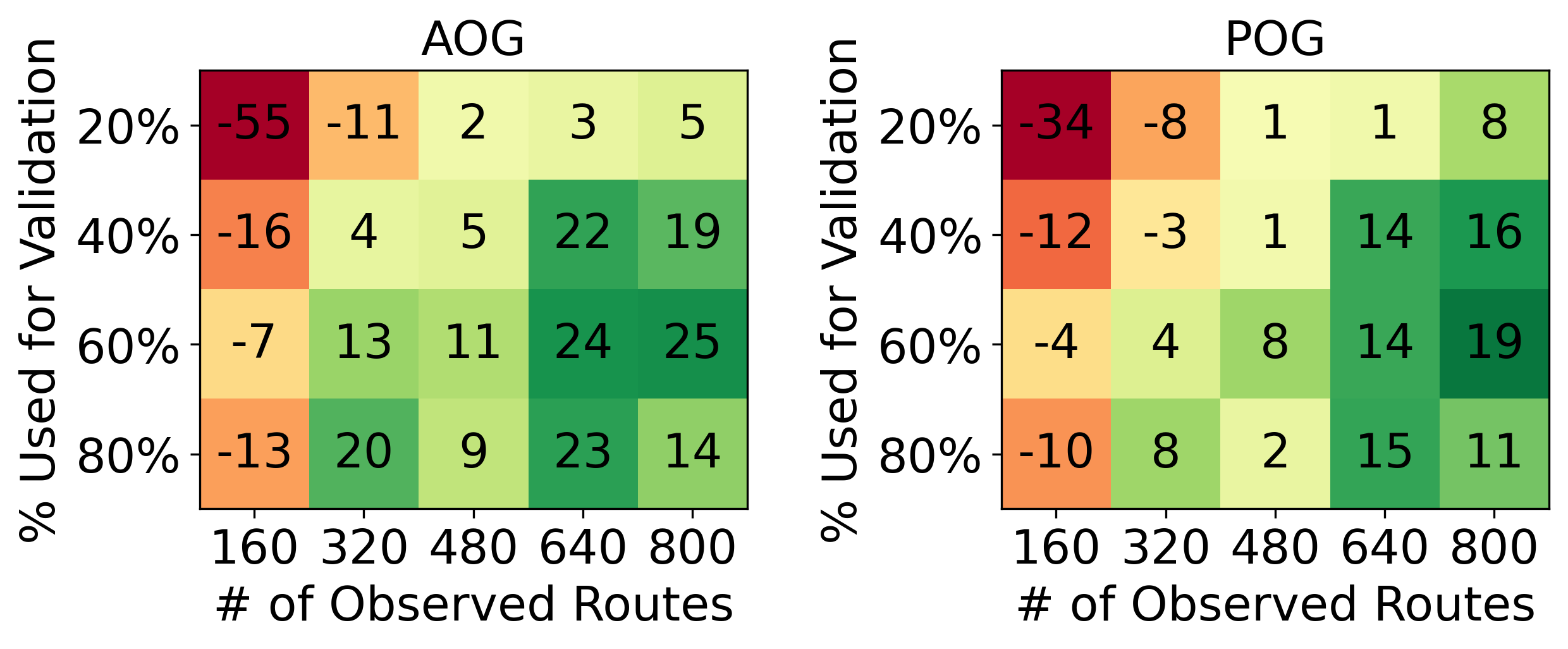}
    \caption{Percentage reduction in AOG and POG when using the conformal IO vs classic IO.}
    \label{fig:size_effect}
\end{figure}

\textbf{Important hyper-parameters.} Finally, we provide empirical evidence that sheds light on the choice of two important hyper-parameters in conformal IO: i) \textit{confidence level $\gamma$}, and ii) \textit{train-validation split ratio}. Regarding $\gamma$, as shown in Figure \ref{fig:pog_aog}, the performance of conformal IO improves quickly as $\gamma$ increases from 0 to 50\% and remains stable and even worsens slightly after that. Hence, it is possible to improve the performance of conformal IO by carefully tuning $\gamma$ using cross-validation. However, this requires an additional validation dataset. If such a dataset is unavailable, setting $\gamma$ to a relatively large value usually yields decent performance, which aligns with our theoretical analysis. Regarding the \textit{train-validation split ratio}, intuitively, both the estimation and calibration can benefit from more data. However, when the dataset is small, we need to strike a balance between these two steps aiming to achieve lower AOG and POG. We implement conformal IO for the shortest path problem under different dataset sizes ($N_\textrm{train} + N_\textrm{val} \in \{160, 320, \ldots, 800\}$) and train-validation split ratios ($N_\textrm{val}/ (N_\textrm{train} + N_\textrm{val})\in \{20\%, 40\%, 60\%, 80\% \}$). As shown in Figure \ref{fig:size_effect}, when the dataset is small (160), there is no benefit of using conformal IO simply because we do not have enough data to obtain a good point estimation and a good uncertainty set at the same time. However, the performance of classic IO quickly plateaus as the dataset grows. When given a mid- or large-sized dataset, we can generally benefit from using more data points in the calibration step, echoing our theoretical analysis.

\section{Conclusion and Future Work} \label{sec:conclusion}

In this paper, we propose conformal IO, a novel approach to recommending high-quality decisions that align with human intuition. We present the first approach to learning uncertainty sets from decision data, which is then utilized in a robust optimization model to prescribe new decisions. We prove that conformal IO achieves bounded optimality gaps, as measured by the ground-truth parameters and the DM's perceived parameters, and demonstrate its strong empirical performance via numerical studies. Finally, we highlight several challenges that underscore future research directions. First, we focus on objective functions that are linear in the unknown parameters. While such objectives are ubiquitous in practice, it is of interest to extend our results for general convex objectives. Second, the robust forward problem, while leading to better decisions, is computationally more expensive than the original forward problem. Future research can be done to accelerate its solution process.


\bibliography{paper.bbl}

\begin{thebibliography}{}

\bibitem[Ahuja and Orlin, 2001]{ahuja2001inverse}
Ahuja, R.~K. and Orlin, J.~B. (2001).
\newblock Inverse optimization.
\newblock {\em Operations Research}, 49(5):771--783.

\bibitem[Aswani et~al., 2018]{aswani2018inverse}
Aswani, A., Shen, Z.-J., and Siddiq, A. (2018).
\newblock Inverse optimization with noisy data.
\newblock {\em Operations Research}, 66(3):870--892.

\bibitem[Babier et~al., 2020]{babier2020knowledge}
Babier, A., Mahmood, R., McNiven, A.~L., Diamant, A., and Chan, T.~C. (2020).
\newblock Knowledge-based automated planning with three-dimensional generative adversarial networks.
\newblock {\em Medical Physics}, 47(2):297--306.

\bibitem[Berthet et~al., 2020]{berthet2020learn}
Berthet, Q., Blondel, M., Teboul, O., Cuturi, M., Vert, J.-P., and Bach, F. (2020).
\newblock Learning with differentiable pertubed optimizers.
\newblock In {\em Advances in Neural Information Processing Systems}, volume~33, pages 9508--9519.

\bibitem[Bodur et~al., 2022]{bodur2022inverse}
Bodur, M., Chan, T.~C., and Zhu, I.~Y. (2022).
\newblock Inverse mixed integer optimization: Polyhedral insights and trust region methods.
\newblock {\em INFORMS Journal on Computing}, 34(3):1471--1488.

\bibitem[Burton et~al., 2020]{burton2020systematic}
Burton, J.~W., Stein, M.-K., and Jensen, T.~B. (2020).
\newblock A systematic review of algorithm aversion in augmented decision making.
\newblock {\em Journal of Behavioral Decision Making}, 33(2):220--239.

\bibitem[Chan et~al., 2023a]{chan2023got}
Chan, T.~C., Mahmood, R., O’Connor, D.~L., Stone, D., Unger, S., Wong, R.~K., and Zhu, I.~Y. (2023a).
\newblock Got (optimal) milk? pooling donations in human milk banks with machine learning and optimization.
\newblock {\em Manufacturing \& Service Operations Management}, 0(0).

\bibitem[Chan et~al., 2014]{chan2014generalized}
Chan, T. C.~Y., Craig, T., Lee, T., and Sharpe, M.~B. (2014).
\newblock Generalized inverse multiobjective optimization with application to cancer therapy.
\newblock {\em Operations Research}, 62(3):680--695.

\bibitem[Chan and Kaw, 2020]{chan2020inverse}
Chan, T. C.~Y. and Kaw, N. (2020).
\newblock Inverse optimization for the recovery of constraint parameters.
\newblock {\em European Journal of Operational Research}, 282(2):415--427.

\bibitem[Chan et~al., 2019]{chan2019inverse}
Chan, T. C.~Y., Lee, T., and Terekhov, D. (2019).
\newblock Inverse optimization: Closed-form solutions, geometry, and goodness of fit.
\newblock {\em Management Science}, 65(3):1115--1135.

\bibitem[Chan et~al., 2023b]{chan2023inverse}
Chan, T. C.~Y., Mahmood, R., and Zhu, I.~Y. (2023b).
\newblock Inverse optimization: Theory and applications.
\newblock {\em Operations Research}, 0(0).

\bibitem[Chen et~al., 2023]{chen2023understanding}
Chen, V., Liao, Q.~V., Wortman~Vaughan, J., and Bansal, G. (2023).
\newblock Understanding the role of human intuition on reliance in human-{AI} decision-making with explanations.
\newblock In {\em Proceedings of the ACM on Human-Computer Interaction}, volume~7, pages 1--32.

\bibitem[Chenreddy et~al., 2022]{chenreddy2022data}
Chenreddy, A.~R., Bandi, N., and Delage, E. (2022).
\newblock Data-driven conditional robust optimization.
\newblock In {\em Advances in Neural Information Processing Systems}, volume~35, pages 9525--9537.

\bibitem[Christiano et~al., 2017]{christiano2017deep}
Christiano, P.~F., Leike, J., Brown, T., Martic, M., Legg, S., and Amodei, D. (2017).
\newblock Deep reinforcement learning from human preferences.
\newblock In {\em Advances in neural information processing systems}, volume~30.

\bibitem[Delage and Ye, 2010]{delage2010distributionally}
Delage, E. and Ye, Y. (2010).
\newblock Distributionally robust optimization under moment uncertainty with application to data-driven problems.
\newblock {\em Operations Research}, 58(3):595--612.

\bibitem[Donahue et~al., 2023]{donahue2023two}
Donahue, K., Kollias, K., and Gollapudi, S. (2023).
\newblock When are two lists better than one?: Benefits and harms in joint decision-making.
\newblock {\em arXiv preprint arXiv:2308.11721}.

\bibitem[Dong et~al., 2018]{dong2018generalized}
Dong, C., Chen, Y., and Zeng, B. (2018).
\newblock Generalized inverse optimization through online learning.
\newblock In {\em Advances in Neural Information Processing Systems}, volume~31.

\bibitem[Dong and Zeng, 2021]{dong2021wasserstein}
Dong, C. and Zeng, B. (2021).
\newblock Wasserstein distributionally robust inverse multiobjective optimization.
\newblock In {\em Proceedings of the AAAI Conference on Artificial Intelligence}, number~7 in 35, pages 5914--5921.

\bibitem[Elmachtoub et~al., 2023]{elmachtoub2023estimate}
Elmachtoub, A.~N., Lam, H., Zhang, H., and Zhao, Y. (2023).
\newblock Estimate-then-optimize versus integrated-estimation-optimization: A stochastic dominance perspective.
\newblock {\em arXiv preprint arXiv:2304.06833}.

\bibitem[Gao and Kleywegt, 2023]{gao2023distributionally}
Gao, R. and Kleywegt, A. (2023).
\newblock Distributionally robust stochastic optimization with wasserstein distance.
\newblock {\em Mathematics of Operations Research}, 48(2):603--655.

\bibitem[Hu et~al., 2022]{hu2022deepeta}
Hu, X., Cirit, O., Binaykiya, T., and Hora, R. (2022).
\newblock {DeepETA}: How uber predicts arrival times using deep learning.
\newblock Uber Engineering Blog.
\newblock Available at \url{https://www.uber.com/en-CA/blog/deepeta-how-uber-predicts-arrival-times/}. Accessed: 2024-01-19.

\bibitem[Ji et~al., 2023]{ji2023ai}
Ji, J., Qiu, T., Chen, B., Zhang, B., Lou, H., Wang, K., Duan, Y., He, Z., Zhou, J., Zhang, Z., et~al. (2023).
\newblock {AI} alignment: A comprehensive survey.
\newblock {\em arXiv preprint arXiv:2310.19852}.

\bibitem[Liu et~al., 2023]{liu2023algorithm}
Liu, M., Tang, X., Xia, S., Zhang, S., Zhu, Y., and Meng, Q. (2023).
\newblock Algorithm aversion: Evidence from ridesharing drivers.
\newblock {\em Management Science}, 0(0).

\bibitem[Mandi et~al., 2022]{mandi2022decision}
Mandi, J., Bucarey, V., Tchomba, M. M.~K., and Guns, T. (2022).
\newblock Decision-focused learning: through the lens of learning to rank.
\newblock In {\em International Conference on Machine Learning}, pages 14935--14947. PMLR.

\bibitem[Merch{\'a}n et~al., 2022]{merchan20222021}
Merch{\'a}n, D., Arora, J., Pachon, J., Konduri, K., Winkenbach, M., Parks, S., and Noszek, J. (2022).
\newblock 2021 {Amazon} last mile routing research challenge: Data set.
\newblock {\em Transportation Science}, 0(0).

\bibitem[Mohajerin~Esfahani et~al., 2018]{mohajerin2018dataIO}
Mohajerin~Esfahani, P., Shafieezadeh-Abadeh, S., Hanasusanto, G.~A., and Kuhn, D. (2018).
\newblock Data-driven inverse optimization with imperfect information.
\newblock {\em Mathematical Programming}, 167:191--234.

\bibitem[Mohri et~al., 2018]{mohri2018foundations}
Mohri, M., Rostamizadeh, A., and Talwalkar, A. (2018).
\newblock {\em Foundations of machine learning}.
\newblock MIT press.

\bibitem[Ng and Russell, 2000]{ng2000algorithms}
Ng, A.~Y. and Russell, S.~J. (2000).
\newblock Algorithms for inverse reinforcement learning.
\newblock In {\em Proceedings of the Seventeenth International Conference on Machine Learning}, page 663–670.

\bibitem[Nguyen, 2015]{nguyen2015uber}
Nguyen, T. (2015).
\newblock {ETA} phone home: How uber engineers an efficient route.
\newblock Uber Engineering Blog.
\newblock Available at \url{https://www.uber.com/en-CA/blog/engineering-routing-engine/}. Accessed: 2024-01-19.

\bibitem[Rafailov et~al., 2024]{rafailov2024direct}
Rafailov, R., Sharma, A., Mitchell, E., Manning, C.~D., Ermon, S., and Finn, C. (2024).
\newblock Direct preference optimization: Your language model is secretly a reward model.
\newblock {\em Advances in Neural Information Processing Systems}, 36.

\bibitem[R{\"o}nnqvist et~al., 2017]{ronnqvist2017calibrated}
R{\"o}nnqvist, M., Svenson, G., Flisberg, P., and J{\"o}nsson, L.-E. (2017).
\newblock Calibrated route finder: Improving the safety, environmental consciousness, and cost effectiveness of truck routing in sweden.
\newblock {\em Interfaces}, 47(5):372--395.

\bibitem[Shafer and Vovk, 2008]{shafer2008tutorial}
Shafer, G. and Vovk, V. (2008).
\newblock A tutorial on conformal prediction.
\newblock {\em Journal of Machine Learning Research}, 9(3).

\bibitem[Sun et~al., 2023]{sun2023predict}
Sun, C., Liu, L., and Li, X. (2023).
\newblock Predict-then-calibrate: A new perspective of robust contextual {LP}.
\newblock In {\em Advances in Neural Information Processing Systems}.

\bibitem[Sun et~al., 2022]{sun2022predicting}
Sun, J., Zhang, D.~J., Hu, H., and Van~Mieghem, J.~A. (2022).
\newblock Predicting human discretion to adjust algorithmic prescription: A large-scale field experiment in warehouse operations.
\newblock {\em Management Science}, 68(2):846--865.

\bibitem[Tan et~al., 2020]{tan2020learn}
Tan, Y., Terekhov, D., and Delong, A. (2020).
\newblock Learning linear programs from optimal decisions.
\newblock In {\em Advances in Neural Information Processing Systems}, volume~33, pages 19738--19749.

\bibitem[Tang and Khalil, 2022]{tang2022pyepo}
Tang, B. and Khalil, E.~B. (2022).
\newblock Pyepo: A pytorch-based end-to-end predict-then-optimize library for linear and integer programming.
\newblock {\em arXiv preprint arXiv:2206.14234}.

\bibitem[Vovk et~al., 2005]{vovk2005algorithmic}
Vovk, V., Gammerman, A., and Shafer, G. (2005).
\newblock {\em Algorithmic learning in a random world}, volume~29.
\newblock Springer.

\bibitem[Wainwright, 2019]{wainwright2019high}
Wainwright, M.~J. (2019).
\newblock {\em High-dimensional statistics: A non-asymptotic viewpoint}, volume~48.
\newblock Cambridge University Press.

\bibitem[Wilder et~al., 2019]{wilder2019melding}
Wilder, B., Dilkina, B., and Tambe, M. (2019).
\newblock Melding the data-decisions pipeline: Decision-focused learning for combinatorial optimization.
\newblock In {\em Proceedings of the AAAI Conference on Artificial Intelligence}, volume~33, pages 1658--1665.

\bibitem[Wirth et~al., 2017]{wirth2017survey}
Wirth, C., Akrour, R., Neumann, G., and F{\"u}rnkranz, J. (2017).
\newblock A survey of preference-based reinforcement learning methods.
\newblock {\em Journal of Machine Learning Research}, 18(136):1--46.

\bibitem[Wu et~al., 2024]{wu2024inverse}
Wu, F., Ke, J., and Wu, A. (2024).
\newblock Inverse reinforcement learning with the average reward criterion.
\newblock In {\em Advances in Neural Information Processing Systems}, volume~36.

\end{thebibliography}
\bibliographystyle{apalike}

\newpage
\appendix
\onecolumn

\section{Omitted Statements and Proofs in Section \ref{sec:prelim}} \label{sec:proof_prelim}

\subsection{Poof of Lemma \ref{lemma:2d_exp_asy_est}}


\proof{} We first show that $\hat{\bfx} \in \{(0, 1), (u, 0)\}$ almost surely. Let $\delta_u := \arccos\left(1 / \sqrt{1 + u^2}\right)$, so $\cos\delta_u = 1 / \sqrt{1 + u^2}$ and $\sin\delta_u = u/ \sqrt{1 + u^2}$. It is easy to verify that, when $\hat{\bstheta}_k \in \bsTheta_1 := \{(\cos\delta, \sin\delta) \,|\, \delta\in (0, \delta_u]\}$, we have $\hat{\bfx}_k = \tilde{\bfx}(\hat{\bstheta}_k, u) = (0, 1)$ almost surely; When $\hat{\bstheta}_k \in \bsTheta_2 := \{(\cos\delta, \sin\delta) \,|\, \delta\in (\delta_u, \pi/2)\}$, we have $\hat{\bfx}_k = \tilde{\bfx}(\hat{\bstheta}_k, u) = (u, 0)$ almost surely. Since $\hat{\theta}_k$ is uniformly distributed in $\hat{\bstheta} \in \bsTheta = \bsTheta_1 \cup \bsTheta_2$, the distribution of $\hat{\bfx}_k$ is
\begin{equation} \label{prob:saa_io_sub}
    \hat{\bfx}_k = 
    \begin{cases}
        & (0, 1), \textrm{ w.p. } 2\delta_u/\pi \\
        & (u, 0), \textrm{ w.p. } (\pi - 2\delta_u)/\pi.
    \end{cases}
\end{equation}

Given a sample set $\mD=\left\{\bfu_k, \hat{\bfx}_k\right\}_{k\in [N]}$, let $N_1$ and $N_2$, respectively, denote the numbers of $(0, 1)$ and $(u, 0)$ in $\mD$. We next show that when $N_1 > 0$ and $N_2 > 0$, $\theta_u$ is the unique optimal solution to $\mathbf{IO}(\mD)$. Specifically, in Example \ref{ex:2d}, $\mathbf{IO}(\mD)$ is presented as follows. 
\begin{equation}
    \bar{\bstheta}_N :=
    \argmin_{\bstheta \in \bsTheta} \quad
        \frac{N_1}{N} l_1(\bstheta)
        + \frac{N_2}{N} l_2(\bstheta)
\end{equation}
where
\begin{equation}
    l_1(\bstheta) = 
    \begin{cases}
        0, 
            & \textrm{ if } \bstheta\in \bsTheta_1, \\
        \theta_2 - u\theta_1, 
            & \textrm{ if } \bstheta\in \bsTheta_2,
    \end{cases}
\end{equation}
and\begin{equation}
    l_2(\bstheta) = 
    \begin{cases}
        u\theta_1 - \theta_2, 
            & \textrm{ if } \bstheta\in \bsTheta_1, \\
        0, 
            & \textrm{ if } \bstheta\in \bsTheta_2.
    \end{cases}
\end{equation}

A simple calculation gives that when $N_1 > 0$ and $N_2 > 0$, the minimum is 0 which occurs uniquely at $\bstheta = (\cos\delta_u, \sin\delta_u)$; When $N_2 = 0$, the minimum is 0 which occurs when $\bstheta \in \bsTheta_1$; When $N_1 = 0$, the minimum is 0 which occurs when $\bstheta \in \bsTheta_2$. Therefore, we have
\begin{equation}
    \bbP(N_1 N_2 > 0) 
    \leq \bbP\left(\bar{\bstheta}_N = (\cos\delta_u, \sin\delta_u)\right)
    \leq 1.
\end{equation}

Given the probability distribution given in Equation~\eqref{prob:saa_io_sub} and that $\mD$ is generated using i.i.d. samples from $\bbP_\theta$, we have
\begin{equation}
    \bbP(N_1 N_2 > 0) = 1 
    - \left(\frac{2\delta_u}{\pi}\right)^N
    - \left(1 - \frac{2\delta_u}{\pi}\right)^N,
\end{equation}
which converges to 1 as $N$ goes to infinity. Therefore, we conclude that $\bbP(\bar{\bstheta}_N = (\cos\delta_u, \sin\delta_u))$ converges to 1 as $N$ goes to infinity.
\endproof

\subsection{Proof of Proposition \ref{prop:unbd_pog_aog}}


\proof{} We first consider the AOG and POG of the decision policy  $\bar{\bfx}_\textrm{IO}(u) := \tilde{\bfx}(\bstheta_u, u)$ separately in the following two lemmas. 

\begin{lemma} \label{lemma:2d_exp_aog}
    In Example \ref{ex:2d}, let $\bar{\bfx}_\mathrm{IO}(u) = \tilde{\bfx}(\bstheta_u, u)$. For any $v\in \bbR_+$ there exists some $\bar{u} > 1$ such that $\mathrm{AOG}(\bar{\bfx}_\mathrm{IO}) > v$ for any $u > \bar{u}_\textrm{AOG}$.
\end{lemma}

\proof{} According to the definition of $\tilde{\bfx}$, we know that $\bar{\bfx}_\textrm{IO}(u)$ is uniformly drawn from
\begin{equation}
    \mX^\textrm{OPT}(\bstheta_u, u) 
    = 
    \left\{
        \left(\frac{ut}{\sqrt{u^2 + 1}}, 1 - \frac{t}{\sqrt{u^2 + 1}} \right)
    \,\middle|\,
        t\in \left[0, \sqrt{u^2 + 1}\right]
    \right\}.
\end{equation}

Since the ground-truth $\bstheta^* = (\cos(\pi/4), \sin(\pi/4))$, the true optimal solution is $\bfx^* = (0, 1)$ with $\tilde{f}(\bstheta^*, u) = \sqrt{2}/2$. Hence, we have
\begin{align}
    \textrm{AOG}(\bar{\bfx}_\textrm{IO}) 
    = \int_0^{\sqrt{u^2 + 1}}
        \frac{\sqrt{2}}{2\sqrt{u^2 + 1}}
        \left(
            1 
            - \frac{t}{\sqrt{u^2 + 1}}
            + \frac{ut}{\sqrt{u^2 + 1}}
        \right)
        dt
        - \frac{\sqrt{2}}{2} 
    = \frac{\sqrt{2}(u - 1)}{4}
\end{align}
Therefore, for any $v \in \bbR_+$, there exists $\bar{u}_\textrm{AOG} = 2\sqrt{2} v + 1$ such that $\mathrm{AOG}(\bar{\bfx}_\textrm{IO}) > v$ for any $u > \bar{u}_\textrm{AOG}$.

\begin{lemma} \label{lemma:2d_exp_pog}
    In Example \ref{ex:2d}, let $\bar{\bfx}_\mathrm{IO}(u) = \tilde{\bfx}(\bstheta_u, u)$. for any $v\in \bbR_+$ there exists some $\bar{u}_\textrm{POG} > 1$ such that $\mathrm{POG}(\bar{\bfx}_\mathrm{IO}) > v$ for any $u > \bar{u}_\textrm{POG}$.
\end{lemma}

\proof{} According to the definition of $\tilde{\bfx}$, $\bar{\bfx}_\textrm{IO}(u)$ is uniformly drawn from
\begin{equation}
    \mX^\textrm{OPT}(\bstheta_u, u) 
    = 
    \left\{
        \left(\frac{ut}{\sqrt{u^2 + 1}}, 1 - \frac{t}{\sqrt{u^2 + 1}} \right)
    \,\middle|\,
        t\in \left[0, \sqrt{u^2 + 1}\right]
    \right\}.
\end{equation}

It is easy to verify that, when $\hat{\bstheta} \in \bsTheta_1 := \{(\cos\delta, \sin\delta) \,|\, \delta\in (0, \delta_u]\}$, we have $\hat{\bfx}_k = \tilde{\bfx}(\hat{\bstheta}, u) = (0, 1)$ with $\tilde{f}(\hat{\bstheta}, u) = \hat{\theta}_2$ almost surely; When $\hat{\bstheta} \in \bsTheta_2 := \{(\cos\delta, \sin\delta) \,|\, \delta\in (\delta_u, \pi/2\}$, we have $\hat{\bfx}_k = \tilde{\bfx}(\hat{\bstheta}, u) = (u, 0)$ with $\tilde{f}(\hat{\bstheta}, u) = u \hat{\theta}_1$ almost surely. Since the optimal solution drawn from $\mX^\textrm{OPT}(\bstheta_u, u)$ is independent of the DM's perception $\hat{\bstheta}$, we have
\begin{align*}
    \mathrm{POG}(\bar{\bfx}_\textrm{IO}) 
    = \  
    & \int_0^{\delta_u} \int_0^{\sqrt{u^2 + 1}}
        \frac{1}{\sqrt{u^2 + 1}}
        \left[
            \frac{ut}{\sqrt{u^2 + 1}} \cos\delta 
            + \left(
                1 - \frac{t}{\sqrt{u^2 + 1}}
            \right)
            \sin\delta
            - \sin\delta 
        \right] 
        dt\, d\delta \\
    & + \int_{\delta_u}^{\pi/2} \int_0^{\sqrt{u^2 + 1}}
        \frac{1}{\sqrt{u^2 + 1}}
        \left[
            \frac{ut}{\sqrt{u^2 + 1}} \cos\delta 
            + \left(
                1 - \frac{t}{\sqrt{u^2 + 1}}
            \right)
            \sin\delta
            -u\cos\delta
        \right] 
        dt\, d\delta \\
    = \  
    & \frac{1}{2}
        \int_0^{\delta_u} \left(u\cos\delta - \sin\delta \right) d\delta
        + \frac{1}{2}
        \int_{\delta_u}^{\pi/2} \left(-u\cos\delta + \sin\delta\right) d\delta \\
    = \  
    & \sqrt{1 + u^2}
        - \frac{u + 1}{2}.\\
    > \  
    & \frac{u-1}{2}
\end{align*}
The inequality holds because $\sqrt{1 + u^2} > u$. Therefore, we have, for any $v \in \bbR_+$, there exists $\bar{u}_\textrm{POG} = 2v + 1$ such that $\mathrm{POG}(\bar{\bfx}_\textrm{IO}) > v$ for any $u > \bar{u}_\textrm{POG}$.
\endproof

Based on Lemmas \ref{lemma:2d_exp_aog} and \ref{lemma:2d_exp_pog}, we conclude that for any $v \in \bbR_+$, there exists $\bar{u} = \max\{\bar{u}_\textrm{AOG}, \bar{u}_\textrm{POG}\} = 2\sqrt{2} v + 1$ such that $\mathrm{AOG}(\bar{\bfx}_\textrm{IO}) > v$ and  $\mathrm{POG}(\bar{\bfx}_\textrm{IO}) > v$ for any $u > \bar{u}$.
\endproof

\subsection{Robustifying the Inverse Problem} \label{sec:alter_robust}

An alternative approach to robustify the classic IO pipeline is to solve a distributionally robust inverse optimization problem. Specifically, consider the following loss function proposed by \citet{mohajerin2018dataIO}.
\begin{definition}
    The distributionally robust sub-optimality loss of $\bstheta$ is given by
    \begin{equation} \label{eq:dro_io_loss}
        \ell_\textrm{DR-S}
        \left( \bstheta \right)
        := 
        \sup_{\bbQ \in \mathfrak{B}^p_r(\hat{\bbP}_{\bfu, \hat{\bfx}})}
        \rho^\bbQ
        \left[
            \ell_\textrm{S} 
            \left(
            \hat{\bfx}, \mX^\textrm{OPT}(\bstheta, \bfu)
            \right)
        \right]
    \end{equation}
    where $\hat{\bbP}_{\bfu, \hat{\bfx}}$ is the sample distribution of $\mD$, $\mathfrak{B}^p_r(\hat{\bbP}_{\bfu, \hat{\bfx}})$ is a $p$-Wasserstain ball of radius $r$ centered at $\hat{\bbP}_{\bfu, \hat{\bfx}}$, and $\rho^\bbQ$ is a risk measure, e.g., the value at risk.
\end{definition}

The \textit{distributionally robust inverse optimization} problem is
\begin{equation}
    \mathbf{DRIO(\mD)}: 
    \underset{\bstheta\in \bsTheta}{\textrm{minimize}}\quad 
        \ell_\textrm{DR-S}(\bstheta).
\end{equation}

As shown by \citet{mohajerin2018dataIO}, the estimated parameters from $\mathbf{DRIO}$ achieve bounded out-of-sample sub-optimality loss with a high probability. However, this does not imply bounded AOG and POG for the decision policy. 

\begin{lemma} \label{lemma:2d_exp_drio_asy}
    In Example \ref{ex:2d}, $\bstheta_u$ is an optimal solution to $\mathbf{DRIO}(\mD)$.
\end{lemma}

\proof{}
As shown in the proof of Lemma \ref{lemma:2d_robust_aog}, when $\alpha\in (0, \pi/2)$, the $\mathbf{RFO}\left( \mC(\bstheta_u, \alpha), u \right)$ has a unique optimal solution $(0, 1)$. So $\bar{\bfx}_\textrm{CIO}(u) = (0, 1)$ almost surely, when $\bar{\bstheta}_N = \bstheta_u$ and $\alpha \in (0, \pi/2)$. It is easy to verify that, when $\hat{\bstheta} \in \bsTheta_1 := \{(\cos\delta, \sin\delta) \,|\, \delta\in (0, \delta_u]\}$, we have $\hat{\bfx}_k = \tilde{\bfx}(\hat{\bstheta}, u) = (0, 1)$ almost surely; When $\hat{\bstheta} \in \bsTheta_2 := \{(\cos\delta, \sin\delta) \,|\, \delta\in (\delta_u, \pi/2)\}$, we have $\hat{\bfx}_k = \tilde{\bfx}(\hat{\bstheta}, u) = (u, 0)$ almost surely. Hence, we have
\begin{equation}
    \mathrm{POG}(\bar{\bfx}_\textrm{CIO}) 
    = \int_{0}^{\delta_u} \frac{\pi}{2} \times 0\ d\delta
        + \int_{\delta_u}^{\pi/2} \frac{\pi}{2} \times \sin\delta\ d\delta
    = - \frac{\pi}{2}\cos\delta \Big|^{\pi/2}_{\delta_u}
    = \frac{\pi}{2\sqrt{1 + u^2}}
    < \frac{\pi}{2\sqrt{2}}.
\end{equation}

The inequality holds because $u > 1$.  

According to Lemma \ref{lemma:2d_exp_asy_est}, we know that $\bbP(\bar{\bstheta}_N = \bstheta_u) \rightarrow 1$ as $N \rightarrow \infty$. So we conclude that, when $\alpha\in (0, \pi/2)$, we have $\bbP\left[\mathrm{POG}(\bar{\bfx}_\textrm{CIO}) < \pi/2\sqrt{2} \right] \rightarrow 1$ as $N \rightarrow \infty$.
\endproof

Lemma \ref{lemma:2d_exp_drio_asy} shows that, in Example \ref{ex:2d}, the estimated parameter from $\mathbf{DRIO}(\mD)$ may still be $\bstheta_u$. Hence, the decision policy is identical to $\bar{\bfx}_\textrm{IO}$ whose AOG and POG can be unbounded. The fundamental reason behind these negative results is the misalignment between the sub-optimality loss and the evaluation metrics. Achieving a low sub-optimality loss means that the suggested and observed decisions are of similar quality as evaluated using the estimated parameters. However, this does not speak to the similarity between these two decisions with respect to the DM's perceived parameters (POG) or the ground-truth parameters (AOG). Therefore, the out-of-sample guarantees on the sub-optimality loss do not translate into bounded AOG or POG.

\subsection{Proof of Lemma \ref{lemma:2d_robust_aog}}

We consider the AOG and POG of conformal IO separately in the following two lemmas.
\begin{lemma}
\label{lemma:2d_robust_aog}
    In Example \ref{ex:2d}, let $\bar{\bfx}_\mathrm{CIO}(u)$ be an optimal solution to $\mathbf{RFO} \left(\mC(\bar{\bstheta}_N, \alpha), u \right)$ where $\bar{\bstheta}_N$ is an optimal solution to $\mathbf{IO}(\mD)$ with the sub-optimality loss~\eqref{eq:sub_opt_loss}. When $\alpha \in (0, \pi/2)$, we have $\bbP\left[\mathrm{AOG}(\bfx_\mathrm{CIO}) = 0\right] \rightarrow 1$ as $N \rightarrow \infty$.
\end{lemma}

\proof{}
We first show that when $\bar{\bstheta} = \bstheta_u$ and $\alpha \in (0, \pi/2)$, $\mathbf{RFO}\left(\mC(\bar{\bstheta}, \alpha), u\right)$ has a unique optimal solution $(0, 1)$. Let $\bfx_1 = (0, 1)$, $\bfx_2 = (0, 2)$, $\bfx_3 = (u, 0)$ and $\bfx_4 = (u, 2)$ denote the four extreme points of the feasible region $\mX(u)$, respectively, and
\begin{equation}
    R(\bfx) := \underset{\bstheta \in \mC(\bstheta_u, \alpha)}{ \max} \ \theta_1 x_1 + \theta_2 x_2.
\end{equation}

Since $\mathbf{FO}$ is a linear program, it suffices to show that, when $\alpha \in (0, \pi/2)$, $R(\bfx_1) < \min\{R(\bfx_2), R(\bfx_3), R(\bfx_4)\}$ because, if there exists an optimal solution that is not an extreme point, then there must exist another extreme point $\bfx_i$ such that $R(\bfx_1) = R(\bfx_i)$ where $i \neq 1$. Next, we compare $R(\bfx_1)$ with $R(\bfx_2)$, $R(\bfx_3)$, and $R(\bfx_4)$. 
 
It is easy to verify that 
\begin{equation}
    R(\bfx_1) = 
    \begin{cases}
        \sin(\delta_u + \alpha), 
            & \textrm{ if } \alpha \in (0, \pi/2 - \delta_u], \\
        1, 
            & \textrm{ if } \alpha \in (\pi/2 - \delta_u, \pi/2).
    \end{cases}
\end{equation}

For $\bfx_2$, we have 
\begin{equation}
    R(\bfx_2) = 
    \begin{cases}
        2 \sin(\delta_u + \alpha), 
            & \textrm{ if } \alpha \in (0, \pi/2 - \delta_u], \\
        2, 
            &\textrm{ if } \alpha \in (\pi/2 - \delta_u, \pi/2).
    \end{cases}
\end{equation}
Hence, we have $R(\bfx_1) < R(\bfx_2)$ when $\alpha\in (0, \pi/2)$. 

For $\bfx_3$, we have 
\begin{equation}
    R(\bfx_3) = 
    \begin{cases}
        u\cos(\delta_u - \alpha), 
            & \textrm{ if } \alpha \in (0, \delta_u], \\
        u, 
            & \textrm{ if } \alpha \in (\delta_u, \pi/2).
    \end{cases}
\end{equation}


Since $u > 1$, we have $\pi/2 - \delta_u < \pi/4 <  \delta_u < \pi/2$. We will show that $R(\bfx_1) < R(\bfx_3)$ when $\alpha$ is in $(0, \pi/2 - \delta_u)$, $[\pi/2 - \delta_u, \delta_u)$, and $[\delta_u, \pi/2)$. When $\alpha \in (0, \pi/2 - \delta_u)$, we have
\begin{align*}
    R(\bfx_1) 
    = \ 
    & \sin(\delta_u + \alpha) \\
    = \ 
    & \sin\delta_u \cos\alpha + \cos\delta_u \sin \alpha
    \\ 
    = \ 
    & \frac{u}{\sqrt{1 + u^2}} \cos\alpha
        + \frac{1}{\sqrt{1 + u^2}} \sin\alpha \\
    < \ 
    & \frac{u}{\sqrt{1 + u^2}} \cos\alpha
        + \frac{u^2}{\sqrt{1 + u^2}} \sin\alpha \\
    = \ 
    & u
    \left(
        \frac{1}{\sqrt{1 + u^2}} \cos\alpha
        + \frac{u}{\sqrt{1 + u^2}} \sin\alpha
    \right) \\
    = \ 
    & u \left(\cos\delta_u \cos\alpha + \sin\delta_u \sin\alpha\right) \\
    = \ 
    & u \cos(\delta_u - \alpha) \\
    = \ & R(\bfx_3).
\end{align*}

The second line holds due to the sum of angles identity. The third line holds due to the definition of $\delta_u$. The fourth line holds because $u > 1$. The fifth line is obtained by simple manipulation. The sixth line holds due to the definition of $\delta_u$. The seventh line holds due to the sum of angles identity. 

When $\alpha \in [\pi/2 - \delta_u, \delta_u)$, we have
\begin{align*}
    R(\bfx_1) 
    = \ 
    & 1 \\
    < \ 
    & 1 + \frac{u-1}{u^2 + 1} \\
    = \ 
    & \frac{u}{\sqrt{u^2 + 1}}
    \frac{1}{\sqrt{u^2 + 1}}
    + 
    \frac{u^2}{\sqrt{u^2 + 1}}
    \frac{1}{\sqrt{u^2 + 1}} \\
    = \ 
    & u\cos\delta_u
    \frac{1}{\sqrt{u^2 + 1}}
    + 
    u\sin\delta_u
    \frac{1}{\sqrt{u^2 + 1}} \\
    < \ 
    & u\cos\delta_u\cos\alpha
    + u\sin\delta_u \sin\alpha \\
    = \ 
    & u\cos(\delta_u -\alpha)\\
    = \ 
    & R(\bfx_3).
\end{align*}

The second line holds because $u > 1$. The third line is obtained through simple manipulation. The forth line holds due to the definition of $\delta_u$. For the fifth line, we know that $\alpha \in [\pi/2 - \delta_u, \delta_u) \subseteq [\pi/4 - \pi/2]$ where $\cos\alpha$ is strictly decreasing in $\alpha$ and where $\sin\alpha$ is strictly increasing in $\alpha$. Therefore, $\cos\alpha < \cos\delta_u = 1/\sqrt{u^2 + 1}$ and $\sin\alpha \leq \sin(\pi/2 - \delta_u) = \cos\delta_u = 1/\sqrt{u^2 + 1}$. Hence, the fifth line holds. The sixth line holds due to the sum of angles identity.

When $\alpha \in [\delta_u, \pi/2)$, we have $R(\bfx_1) = 1 < u = R(\bfx_3)$. 

Hence, $R(\bfx_1) < R(\bfx_3)$ when $\alpha\in (0, \pi/2)$. 

For $\bfx_4$, we have 
\begin{equation}
    R(\bfx_4) = \max_{\delta \in \mC(\delta_u, \alpha)}
    u \cos\delta + 2\sin\delta.
\end{equation}

Let $\delta^*_1$ denote the optimal solution to the maximization problem for calculating $R(\bfx_1)$. It is easy to verify that $\delta^*_1 \in (0, \pi/2)$ when $\alpha \in (0, \pi/2)$. So $\cos\delta^*_1 >0$ and $\sin\delta^*_1 > 0$. Hence, we have
\begin{equation}
    R(\bfx_4) = \max_{\delta \in \mC(\delta_u, \alpha)}
    u \cos\delta + 2\sin\delta
    \geq 
    u \cos\delta^*_1 + 2\sin\delta^*_1
    > 
    \sin\delta^*_1
    = R(\bfx_1).
\end{equation}
The first inequality holds because $\delta^*_1$ may not be the maximizer of the problem associated with $\bfx_4$. The second inequality holds because $u >1$, $\cos\delta^*_1 > 0$, and $\sin\delta^*_1 > 0$.

Hence, when $\alpha\in (0, \pi/2)$, $\mathbf{RFO}\left( \mC(\delta_u, \alpha), u \right)$ has a unique optimal solution $\bfx_1$, so $\bar{\bfx}_\textrm{CIO}(u) = (0, 1)$ almost surely. Given that $(0, 1)$ is also the optimal solution to $\mathbf{FO}(\delta^*, u)$, we have $\mathrm{AOG}(\bar{\bfx}_\textrm{CIO}) = 0$ when $\bar{\bstheta}_N = \bstheta_u$ and $\alpha\in (0, \pi/2)$. According to Lemma \ref{lemma:2d_exp_asy_est}, we know that $\bbP(\bar{\bstheta}_N = \bstheta_u) \rightarrow 1$ as $N \rightarrow \infty$. So we conclude that, when $\alpha\in (0, \pi/2)$, we have $\bbP\left[\mathrm{AOG}(\bar{\bfx}_\textrm{CIO}) = 0 \right] \rightarrow 1$ as $N \rightarrow \infty$.

\endproof

\begin{lemma}
\label{lemma:2d_robust_pog}
    In Example \ref{ex:2d}, let $\bar{\bfx}_\mathrm{CIO}(u)$ be an optimal solution to $\mathbf{RFO} \left(\mC(\bar{\bstheta}_N, \alpha), u \right)$ where $\bar{\bstheta}_N$ is an optimal solution to $\mathbf{IO}(\mD)$ with the sub-optimality loss~\eqref{eq:sub_opt_loss}.
    When $\alpha \in (0, \pi/2)$, we have $\bbP\left[\mathrm{POG}(\bar{\bfx}_\mathrm{CIO}) < \pi/2\sqrt{2}\right] \rightarrow 1$ as $N\rightarrow \infty$.
\end{lemma}
\proof{}

As shown in the proof of Lemma \ref{lemma:2d_robust_aog}, when $\alpha\in (0, \pi/2)$, the $\mathbf{RFO}\left( \mC(\bstheta_u, \alpha), u \right)$ has a unique optimal solution $(0, 1)$. So $\bar{\bfx}_\textrm{CIO}(u) = (0, 1)$ almost surely, when $\bar{\bstheta}_N = \bstheta_u$ and $\alpha \in (0, \pi/2)$. It is easy to verify that, when $\hat{\bstheta} \in \bsTheta_1 := \{(\cos\delta, \sin\delta) \,|\, \delta\in (0, \delta_u]\}$, we have $\hat{\bfx}_k = \tilde{\bfx}(\hat{\bstheta}, u) = (0, 1)$ almost surely; When $\hat{\bstheta} \in \bsTheta_2 := \{(\cos\delta, \sin\delta) \,|\, \delta\in (\delta_u, \pi/2)\}$, we have $\hat{\bfx}_k = \tilde{\bfx}(\hat{\bstheta}, u) = (u, 0)$ almost surely. Hence, we have
\begin{equation}
    \mathrm{POG}(\bar{\bfx}_\textrm{CIO}) 
    = \int_{0}^{\delta_u} \frac{\pi}{2} \times 0\ d\delta
        + \int_{\delta_u}^{\pi/2} \frac{\pi}{2} \times \sin\delta\ d\delta
    = - \frac{\pi}{2}\cos\delta \Big|^{\pi/2}_{\delta_u}
    = \frac{\pi}{2\sqrt{1 + u^2}}
    < \frac{\pi}{2\sqrt{2}}.
\end{equation}

The inequality holds because $u > 1$.  

According to Lemma \ref{lemma:2d_exp_asy_est}, we know that $\bbP(\bar{\bstheta}_N = \bstheta_u) \rightarrow 1$ as $N \rightarrow \infty$. So we conclude that, when $\alpha\in (0, \pi/2)$, we have $\bbP\left[\mathrm{POG}(\bar{\bfx}_\textrm{CIO}) < \pi/2\sqrt{2} \right] \rightarrow 1$ as $N \rightarrow \infty$.
\endproof
\section{Proof of Statements in Section \ref{sec:cio}}

\subsection{Definitions} \label{subsec:def}

\begin{definition}[Empirical Rademacher Complexity]
    Let $\mF$ be a class of functions mapping from $\mZ = \{Z_1, Z_2, \ldots, Z_m\}$ to $[a, b]$ and $\mD$ be a fixed sample of size $N$ with elements in $\mZ$, then the empirical Rademacher Complexity of $\mF$ with respect to the sample $\mD$ is defined as 
    \begin{equation}
        \hat{\mathfrak{R}}_{\mD}(\mF) := 
        \bbE_{\bssigma}
        \left[
            \sup_{f\in \mF}
            \frac{1}{N} 
            \sum_{i\in [N]} \sigma_i f(Z_i)
        \right]
    \end{equation}
    where $\bssigma = (\sigma_1, \sigma_2, \ldots, \sigma_N)^\intercal$ with $\sigma_i$'s being independent uniform random variables taking values in $\{-1, 1\}$.
\end{definition}

\begin{definition}[Rademacher Complexity] \label{def:rad_comp}
    Let $\bbP$ denote the distribution according to which samples are drawn. For any integer $N \geq 1$, the Rademacher complexity of a function class $\mF$ is the expectation of the empirical Rademacher complexity over the samples of size $N$ drawn from $\bbP$:
    \begin{equation}
        \mathfrak{R}_N(\mF) := 
            \bbE_{\mD \sim \bbP^N}
                \left[\hat{\mathfrak{R}}_{\mD}(\mF) \right]
    \end{equation}
\end{definition}

\begin{definition}[Growth Function] \label{def:growth_func}
    Let $\mH$ be a class of functions that take values in $\{-1, 1\}$. The growth function $\Pi_\mH: \bbN \rightarrow \bbN$ for $\mH$ is defined as 
    \begin{equation}
        \Pi_\mH(N) :=
        \max_{(Z_1, Z_2, \ldots, Z_N) \in \mZ^N}
        \left|
            \left\{
                \left(h(Z_1), h(Z_2), \ldots, h(Z_N)\right)
            \,\middle|\,
                h \in \mH    
            \right\}
        \right|
    \end{equation}
    which measures the maximum number of distinct ways in which $N$ data points in $\mZ$ can be classified using the function class $\mH$.
\end{definition}

\subsection{Useful Lemmas} \label{subsec:lemma}

\begin{lemma}[Corollary 3.1 in \citet{mohri2018foundations}]
    \label{lemma:massert}
    Let $\mH$ be a class of functions taking values in $\{1, -1\}$, then, for any integer $N \geq 1$, the following holds
    \begin{equation}
        \mathfrak{R}_N(\mH) 
        \leq
        \sqrt{\frac{2\log \Pi_\mH(N)}{N}}.
    \end{equation}
\end{lemma}

\begin{lemma}[Theorem 4.10 in \citet{wainwright2019high}]
    \label{lemma:rademacher}
   For any $b$-uniformly bounded class of functions $\mF$, any positive integer $N \geq 1$, and any scalar $\delta \geq 0$, with probability at least $1 - \exp\left(-N\delta^2/ (2b^2)\right)$, we have
    \begin{equation}
        \underset{f\in \mF}{\sup}
        \left|
            \frac{1}{N}\sum_{i\in [N]}f(X_i) 
            - \bbE\left[f(X_i)\right]
        \right|
        \leq 
        2 \mathfrak{R}_N(\mF) + \delta
    \end{equation}
    where $\mathfrak{R}(\mF)$ denotes the Rademacher complexity of the function class $\mF$.
\end{lemma}

\subsection{Proof of Theorem \ref{thm:cali_reform_cone}} 
\proof{} We first present the extensive formulation of Problem~\eqref{prob:calib}. For convenience, we define $\hat{\bsTheta}_k := \bsTheta^\textrm{OPT}(\bfu_k, \hat{\bfx}_k)$ for any $k\in [N]$. When $\alpha \in [0, \pi]$, $\cos\alpha$ is a strictly decreasing in $\alpha$. Therefore, minimizing $\alpha$ is equivalent to maximizing the value of $\cos\alpha$. We can replace the decision variable $\alpha$ in Problem~\eqref{prob:calib} with a new decision variable $c := \cos\alpha$ with an additional constraint $t$ with $-1 \leq c \leq 1$. In addition, we introduce a new set of decision variables  $y_k\in \{0, 1\}$ that indicate if $\hat{\bsTheta}_k$ intersects with the learned uncertainty set ($=1$) or not ($=0$) for any $k\in \mK_\textrm{val}$. Problem~\eqref{prob:calib} can be presented as follows.  
\begin{subequations} \label{prob:cali_ext}
\begin{align}
    \underset{c, \{\bstheta_k\}_{k\in \mK_\textrm{val}}, \{y_k\}_{k\in \mK_\textrm{val}}}{\textrm{maximize}} \quad
    & c \\
    \textrm{subject to} \quad 
    \label{ext_cali:inv_fea}
    & \hat{\bfx}_k \in \mX^\textrm{OPT}(\bstheta_k, \bfu_k),
        \quad \forall k \in \mK_\textrm{val} \\
    \label{ext_cali:intersect}
    & \bstheta_k^\intercal \bar{\bstheta}
        \geq c + 2(y_k - 1), 
        \quad \forall k \in \mK_\textrm{val}
        \\
    \label{ext_cali:threshold}
    & \sum_{k\in \mK_\textrm{val}} y_k 
        \geq \lceil \gamma (N_\textrm{val} + 1) \rceil \\
    \label{ext_cali:regu}
    & \|\bstheta_k\|_2 = 1, 
        \quad \forall k\in \mK_\textrm{val} \\
    \label{ext_cali:t}
    & -1\leq c \leq 1 \\
    \label{ext_cali:y}
    & y_k \in \{0, 1\}, 
        \quad \forall k\in \mK_\textrm{val}.
\end{align}
\end{subequations}
Constraints~\eqref{ext_cali:inv_fea} ensure that $\bstheta_k$ is a member of $\hat{\bsTheta}_k$ for any $k\in \mK_\textrm{val}$. Constraints~\eqref{ext_cali:intersect} decide if $\bstheta_k$ should be taken into account when calculating the maximal cosine value $c$ based on if $\hat{\bsTheta}_k$ intersects with $\mC$. Constraint~\eqref{ext_cali:threshold} ensures that $\mC$ intersects with at least $\lceil \gamma (N_\textrm{val} + 1) \rceil$ inverse feasible sets. Constraint~\eqref{ext_cali:regu} enforces $\bstheta_k$ to be on the unit sphere as defined in Equation~\eqref{def:cone}. Constraints~\eqref{ext_cali:t}--\eqref{ext_cali:y} specify the ranges of the decision variables. 

Observing that the objective of Problem~\eqref{prob:cali_ext} is to maximize $c$ and that decision variables $\bstheta_k$ of different data points only interact in Constraints~\eqref{ext_cali:intersect}. We can re-write Problem~\eqref{prob:cali_ext} as 
\begin{subequations} 
\begin{align} \label{prob:cali_master}
    \textrm{maximize} \quad
    & c \\
    \textrm{subject to} \quad 
    & c \leq c_k - 2(y_k - 1), 
        \quad \forall k \in \mK_\textrm{val} \\
    & \sum_{k\in \mK_\textrm{val}} y_k 
        \geq \lceil \gamma (N_\textrm{val} + 1)\rceil \\
    & -1\leq c \leq 1 \\
    & y_k \in \{0, 1\}, 
        \quad \forall k\in \mK_\textrm{val},
\end{align}
\end{subequations}
where
\begin{subequations} \label{prob:max_cos}
\begin{align}
    c_k :=
    \underset{\bstheta_k}{\mathrm{maximize}} \quad 
    & \bstheta_k^\intercal \bar{\bstheta} \\
    \mathrm{subject}\textrm{ }\mathrm{to} \quad
    & \hat{\bfx}_k \in \mX^\textrm{OPT}(\bstheta_k, \bfu_k) \\
    & \|\bstheta_k\|_2 \leq 1. \label{max_cos:conv_reg}
\end{align}
\end{subequations}

Note that we replace Constraints~\eqref{ext_cali:regu} with Constraints~\eqref{max_cos:conv_reg} because the objective of Problem~\eqref{prob:max_cos} is to maximize the inner product of $\bstheta_k$ and $\bar{\bstheta}$, so the maximum only occurs when $\|\bstheta_k\|_2 = 1$. We further observe that the optimal solution to Problem~\eqref{prob:cali_master} is to set $y_k = 1$ for all $k$ such that $c_k \geq \Gamma_\tau\left(\left\{ c_k\right\}_{k\in \mK_\textrm{val}}\right)$ and $y_k = 0$ otherwise. Therefore, the optimal objective value of Problem~\eqref{prob:cali_master} is $c = \Gamma_\tau\left(\left\{ c_k\right\}_{k\in \mK_\textrm{val}}\right)$ 
corresponding to $\alpha_\gamma = \arccos \Gamma_\tau\left(\left\{ c_k\right\}_{k\in \mK_\textrm{val}}\right)$.
\endproof

\subsection{Proof of Lemma \ref{lemma:lipschiz}}
\proof{Proof.}
For any fixed $\bfx$, we have
\begin{align*}
    f\left(\bstheta, \hat{\bfx}\right) - f\left(\bstheta', \hat{\bfx}\right)
    = 
    & \sum_{i\in [d]} 
        \left(\theta_i - \theta'_i\right) 
        f_i\left(\hat{\bfx}\right)  \\
    \leq 
    &  \sqrt{\sum_{i\in [d]} f_i^2(\hat{\bfx})}
        \sqrt{\sum_{i\in [d]} (\theta_i - \theta'_i)^2} \\
    = 
    & \nu(\hat{\bfx}) \|\bstheta - \bstheta'\|_2
\end{align*}
where $\nu(\hat{\bfx}) := \sqrt{\sum_{i\in [d]} f_i^2(\hat{\bfx})}$. The inequality follows the Cauchy-Schwartz inequality.
\endproof

\subsection{Proof of Theorem \ref{thm:valid}} 

\proof{} We first prove the learned uncertainty set is conservatively valid. Following the conformal prediction language used by \citet{vovk2005algorithmic}, we define a conformality measure of each data point,i.e. an observed decision and exogenous parameter pair, $A_{\bar{\bstheta}}: \bbR^{n} \times \mU \rightarrow \bbR_+$ as follows
\begin{subequations}
\begin{align}
    A_{\bar{\bstheta}}(\hat{\bfx}, \bfu) := 
        \underset{\bstheta}{\mathrm{maximize}} \quad 
        & \bstheta^\intercal \bar{\bstheta} \\
        \mathrm{subject}\textrm{ }\mathrm{to} \quad
        & \bstheta \in \bsTheta^\textrm{OPT}(\hat{\bfx}, \bfu) \\
        & \|\bstheta\|_2 \leq 1.
\end{align}
\end{subequations}
We note that $c_k = A_{\bar{\bstheta}}(\hat{\bfx}_k, \bfu_k)$ for any $k\in \mK_\textrm{val}$ where $c_k$ is defined in Theorem \ref{thm:cali_reform_cone}. Let $\tau = \left\lceil \gamma(N_\textrm{val} + 1) \right\rceil$, and $\mA := \{A_{\bar{\bstheta}}(\hat{\bfx}_k, \bfu_k) \}_{k\in \mK_\mathrm{val}}$, or equivalently, $\mA := \left\{c_k \right\}_{k\in \mK_\textrm{val}}$. Due to the definition of $\mC\left(\bar{\bstheta}, \alpha\right)$ and that $\alpha$ is chosen such that $\cos\alpha = \Gamma^\tau (\mA)$, the event ``$\bsTheta^\textrm{OPT}(\hat{\bfx}, \bfu) \cap \mC(\bar{\bstheta}, \alpha) \neq \emptyset$'' is equivalent to ``$A_{\bar{\bstheta}}(\hat{\bfx}, \bfu) \geq \Gamma^\tau (\mA)$'', so
\begin{equation} \label{eq:prob_equ}
    \bbP\left(
        \bsTheta^\textrm{OPT}(\hat{\bfx}, \bfu) 
        \cap \mC(\bar{\bstheta}, \alpha) 
        \neq \emptyset
    \right)
    =
    \bbP\left(
        A_{\bar{\bstheta}}(\hat{\bfx}, \bfu) 
        \geq \Gamma^\tau (\mA)
    \right).
\end{equation}

Assumption \ref{asm:iid} implies that the dataset $\mD' = \mD_\textrm{val} \cup \{(\hat{\bfx}, \bfu)\}$ is exchangeable, i.e. the ordering of the data points in $\mD'$ does not affect its joint probability distribution \citep{shafer2008tutorial}. Therefore, the rank (from high to low) of $A_{\bar{\bstheta}}(\hat{\bfx}, \bfu)$ in $\mA' := \mA \cup \{A_{\bar{\bstheta}}(\hat{\bfx}, \bfu)\}$ is uniformly distributed in $\{1, 2, \ldots, N_\mathrm{val} + 1\}$. So, we have
\begin{align*}
    \gamma 
    \leq \ 
    & \bbP
    \left\{
        A_{\bar{\bstheta}}(\hat{\bfx}, \bfu) 
        \geq
        \Gamma^\tau(\mA')
    \right\} \\
    = \ 
    & 1 - \bbP
    \left\{
        A_{\bar{\bstheta}}(\hat{\bfx}, \bfu) 
        <
        \Gamma^\tau(\mA')
    \right\} \\
    = \ 
    & 1 - \bbP
    \left\{
        A_{\bar{\bstheta}}(\hat{\bfx}, \bfu) 
        <
        \Gamma^\tau(\mA)
    \right\} \\
    = \ 
    & \bbP
    \left\{
        A_{\bar{\bstheta}}(\hat{\bfx}, \bfu) 
        \geq
        \Gamma^\tau(\mA)
    \right\} \\
    = \ 
    & \bbP
    \left\{
        \bsTheta^\textrm{OPT}(\hat{\bfx}, \bfu) \cap \mC(\bar{\bstheta}, \alpha) \neq \emptyset
    \right\}.
\end{align*}

The first line holds due to the definition of $\tau$. We obtain the second line by taking the complement of the event in the first line (inside the probability). The third line holds because $A_{\bar{\bstheta}}(\hat{\bfx}, \bfu)$ can never be strictly smaller than itself, so any elements in $\mA'$ that are strictly smaller than $A_{\bar{\bstheta}}(\hat{\bfx}, \bfu)$ are in $\mA$. Note that this line holds only when $\tau \leq N_\textrm{val}$, which occurs when $\gamma \leq N_\textrm{val} / (N_\textrm{val} + 1)$, because $\mA$ only has $N_\textrm{val}$ elements. We obtain the third line by taking the complement of the event in the second line (inside the probability). The last line holds due to Equation~\eqref{eq:prob_equ}. We note that all the probabilities are over the joint distribution of $\mD_\textrm{val}$ and the new sample, i.e. $\mD'$. 

We next prove that the learned uncertainty set is asymptotically exact. Let $z_k := (\bfu_k, \hat{\bfx}_k)$, $\mZ:= \{z_k\}_{k\in \mK_\textrm{val}}$. We define a function class
\begin{equation}
    \mH = 
    \left\{
        h(z, \alpha) = 
        \mathds{1}
        \left[ 
            \bsTheta^\textrm{OPT}(\hat{\bfx}, \bfu) \cap \mC(\bar{\bstheta}, \alpha)
        \right]
    \, \middle| \,
        \alpha \in (0, \pi)
    \right\}.
\end{equation}

Let $\Pi_\mH$ denote the growth function of $\mH$ as defined in Definition \ref{def:growth_func}. It is easy to verify that
\begin{equation}
    \Pi_{\mH} (N_\textrm{val}) = N_\textrm{val} + 1
\end{equation}
because the value of $h(z, \alpha)$ is monotonically increasing in $\alpha$ for any fixed $z\in \mZ$, so changing the value of $\alpha$ can only leads to $N_\textrm{val} + 1$ different outcomes for a fixed dataset $\mZ$.

Therefore, according to Lemma \ref{lemma:massert}, we have
\begin{equation}
    \label{eq:massert}
    \mathfrak{R}_{N_\textrm{val}}(\mH) \leq \sqrt{\frac{2\log (N_\textrm{val} + 1)}{N_\textrm{val}}}
\end{equation}
where $\mathfrak{R}_{N_\textrm{val}}(\mH)$ denotes the Rademacher complexity of $\mH$ when sample size is $N_\textrm{val}$, as defined in Definition \ref{def:rad_comp}.

We know that the value of $\alpha$ is chosen such that it is the smallest value that satisfies 
\begin{equation}
    \frac{1}{N_\textrm{val}}
        \sum_{k\in \mK_\textrm{val}}h(z_k, \alpha) 
    = 
    \frac{1}{N_\textrm{val}}
    \sum_{k\in \mK_\textrm{val}}\mathds{1}
        \left[ 
            \bsTheta^\textrm{OPT}(\hat{\bfx}_k, \bfu_k) \cap \mC(\bar{\bstheta}, \alpha)
        \right] 
    = 
    \frac{\left\lceil\gamma (N_\textrm{val} + 1)\right\rceil}{N_\textrm{val}},
\end{equation}
so we have
\begin{equation}
    \gamma 
    \leq
    \frac{1}{N_\textrm{val}}
        \sum_{k\in \mK_\textrm{val}}h(z_k, \alpha)
    \leq \gamma + \frac{2}{N_\textrm{val}}.
\end{equation}

The second inequality holds because
\begin{equation}
    \frac{\left\lceil\gamma (N_\textrm{val} + 1)\right\rceil}{N_\textrm{val}}
    = \frac{
        \lfloor\gamma N_\textrm{val} \rfloor
        + \lceil \gamma N_\textrm{val}
            -\lfloor\gamma N_\textrm{val} \rfloor
            + \gamma
            \rceil
        }{N_\textrm{val}}
    \leq 
    \frac{
        \gamma N_\textrm{val}
        + \lceil \gamma N_\textrm{val}
            -\lfloor\gamma N_\textrm{val} \rfloor
            + \gamma
            \rceil
        }{N_\textrm{val}}
    \leq 
    \gamma + \frac{2}{N_\textrm{val}}.
\end{equation}

Since $\mD_\textrm{val}$ is i.i.d. sampled, for any fixed $\alpha$, $\sum_{k\in \mK_\textrm{val}} h(z_k, \alpha) / N_\textrm{val}$ provides a sample average approximation to $\bbE\left[h(z, \alpha)\right]$, which can be interpreted as $ \bbP
\left(
    \bsTheta^\textrm{OPT}(\hat{\bfx}, \bfu) \cap \mC(\bar{\bstheta}, \alpha)
\right)$ for any new sample $(\hat{\bstheta}, \bfu)$ from $\bbP_{\hat{\bstheta}, \bfu}$ and $\hat{\bfx} = \tilde{\bfx}(\hat{\bstheta}, \bfu)$. 

By applying Lemma \ref{lemma:rademacher}, we have, with probability at least $\delta = 1 - 1/N_\textrm{val}$,  
\begin{align}
    \label{eq:rademacher}
    \left|
        \frac{1}{N_\textrm{val}}\sum_{k\in \mK_\textrm{val}}h(z_k, \alpha) 
        - \bbE\left[h(z, \alpha)\right]
    \right|
    \leq 2 \mathfrak{R}_{N_\textrm{val}}(\mH) + \sqrt{\frac{2\log N_\textrm{val}}{N_\textrm{val}}}.
\end{align}

By combing \eqref{eq:massert}--\eqref{eq:rademacher}, we have, with probability at least $1 - 1/N_\textrm{val}$,

\begin{equation}
    \left|
        \bbP
        \left(
            \bsTheta^\textrm{OPT}(\hat{\bfx}, \bfu) \cap \mC(\bar{\bstheta}, \alpha)
        \right) 
        - \gamma
    \right|
    \leq 
    \sqrt{\frac{8\log(N_\textrm{val} + 1) + 2\log N_\textrm{val}}{N_\textrm{val}}}
    + \frac{2}{N_\textrm{val}}.
\end{equation}
\endproof

\subsection{Proof of Theorem \ref{thm:pog}}

We bound the AOG and POG of conformal IO separately in the following two propositions. 

\begin{proposition}[Conformal IO Achieves Bounded POG] \label{prop:pog}
    Let $\bar{\bfx}_\mathrm{CIO}(\bfu)$ be an optimal solution to $\mathbf{RFO}\left(\mC(\bar{\bstheta}, \alpha_1), \bfu \right)$ for any $\bfu\in \mU$, where $\bar{\bstheta} \in \bbR^d$ and $\alpha_1$ are chosen such that, for a new sample $(\bstheta', \bfu')$ from $\bbP_{(\bstheta, \bfu)}$ and $\bfx' = \tilde{\bfx}(\bstheta', \bfu')$, $\bbP \left( \mC(\bar{\bstheta}, \alpha_1) \cap \bsTheta^\mathrm{OPT}(\bfu', \bfx') \neq \emptyset \right) = 1$. If Assumptions \ref{asm:bounded_inv_feas}--\ref{asm:bnd_div} hold, then  
    \begin{equation}
        \mathrm{POG}(\bar{\bfx}_\mathrm{CIO}) 
        \leq 
        (\eta - 2\cos2\alpha_1 + 2)\mu
        + \eta \mu_\mathrm{CIO}
    \end{equation}
    where $\mu := \bbE [\nu(\tilde{\bfx}(\hat{\bstheta}, \bfu))]$ and $\mu_\mathrm{CIO} := \bbE (\nu[ \bar{\bfx}_\mathrm{CIO}(\bfu) ])$.
\end{proposition}

\proof{}
We first bound the perceived optimality gap of a sampled DM. Let $(\hat{\bstheta}, \bfu )$ be a sample from $\bbP_{(\bstheta, \bfu)}$, $\hat{\bfx} = \tilde{\bfx}(\hat{\bstheta}, \bfu)$, 
$\hat{\bstheta}_\textrm{CIO}(\bfu)$ denote the optimal solution to the inner maximization problem in $\mathbf{RFO}\left(\mC(\bar{\bstheta}, \alpha_1), \bfu\right)$ when the outer decision variable is set to $\hat{\bfx}$,
$\bar{\bstheta}_\textrm{CIO}(\bfu)$ denote the optimal solution to the inner maximization problem in $\mathbf{RFO}\left(\mC(\bar{\bstheta}, \alpha_1), \bfu\right)$ when the outer decision variable is set to $\bar{\bfx}_\textrm{CIO}(\bfu)$,
If $\bsTheta^\mathrm{OPT}\left(\hat{\bfx}, \bfu\right) \cap \mC(\bar{\bstheta}, \alpha_1) \neq \emptyset$, let  $\tilde{\bstheta}$ be an element of $\bsTheta^\mathrm{OPT}\left(\hat{\bfx}, \bfu\right) \cap \mC(\bar{\bstheta}, \alpha_1)$, we have

\begin{align*}
    & f\left(
        \hat{\bstheta}, \bar{\bfx}_\mathrm{CIO}(\bfu)
    \right) 
    - f\left(
        \hat{\bstheta}, \hat{\bfx}
    \right) \\
    \leq \ 
    & f
    \left(
        \tilde{\bstheta}, \bar{\bfx}_\mathrm{CIO}(\bfu)
    \right) 
    - f\left(
        \tilde{\bstheta}, \hat{\bfx}
    \right) 
    + \left[
        \nu(\hat{\bfx})
        + \nu\left(\bar{\bfx}_\mathrm{CIO}(\bfu)\right)  
        \right]
        \left\|\hat{\bstheta} - \tilde{\bstheta} \right\|_2 \\
    \leq \ 
    & f
    \left(
        \tilde{\bstheta}, \bar{\bfx}_\mathrm{CIO}(\bfu)
    \right) 
    - f\left(
        \tilde{\bstheta}, \hat{\bfx}
    \right) 
    + \eta \left[
        \nu(\hat{\bfx})
        + \nu\left(\bar{\bfx}_\mathrm{CIO}(\bfu)\right)  
        \right]  \\
    \leq \ 
    & f\left(\bar{\bstheta}_\mathrm{CIO}(\bfu), \bar{\bfx}_\mathrm{CIO}(\bfu)\right) 
        - f\left(\tilde{\bstheta}, \hat{\bfx}\right) 
        + \eta 
            \left[
                \nu(\hat{\bfx})
                + \nu\left(\bar{\bfx}_\mathrm{CIO}(\bfu)\right)  
            \right] \\
    \leq \ 
    & f\left(\hat{\bstheta}_\mathrm{CIO}(\bfu), \hat{\bfx}\right) 
        - f\left(\tilde{\bstheta}, \hat{\bfx}\right) 
        + \eta 
            \left[
                \nu(\hat{\bfx})
                + \nu\left(\bar{\bfx}_\mathrm{CIO}(\bfu)\right)  
            \right] 
            \\            
    \leq \ 
    & \nu(\hat{\bfx}) 
        \left\|
            \hat{\bstheta}_\mathrm{CIO}(\bfu) - \tilde{\bstheta} 
        \right\|_2
        + \eta 
            \left[
                \nu(\hat{\bfx})
                + \nu\left(\bar{\bfx}_\mathrm{CIO}(\bfu)\right)  
            \right] \\
    \leq \ 
    & 2\nu(\hat{\bfx}) (1 - \cos2\alpha_1)+ \eta \left(\nu(\hat{\bfx}) +         \nu\left[\bar{\bfx}_\mathrm{CIO}(\bfu)\right]\right)\\
    = \ 
    & \nu(\hat{\bfx})(\eta - 2\cos2\alpha_1 + 2)
        + \eta \nu\left[\bar{\bfx}_\mathrm{CIO}(\bfu)\right].
\end{align*}

The first line holds due to Lemma \ref{lemma:lipschiz}. 
The second line holds due to assumption \ref{asm:bounded_inv_feas}. 
The third line holds due to the definition of $\bar{\bstheta}_\textrm{CIO}(u)$. The fourth line holds because $\left( \bar{\bfx}_\mathrm{CIO}(\bfu), \bar{\bstheta}_\textrm{CIO}(\bfu) \right)$ is an optimal solution to $\mathbf{RFO}\left(\mC(\bar{\bstheta}, \alpha_1), \bfu\right)$.
The fifth line holds due to Lemma \ref{lemma:lipschiz}. 
The sixth line holds because both $\hat{\bstheta}_\textrm{CIO}(\bfu)$ and $\tilde{\bstheta}$ are in $\mC(\bar{\bstheta}, \alpha_1)$ so the angle between them is no larger than $2\alpha_1$. Since both $\hat{\bstheta}_\textrm{CIO}(\bfu)$ and $\tilde{\bstheta}$ are on the unit sphere, the $L_2$ distance between them are bounded by $2(1 - \cos2\alpha_1)$.

Since $\alpha_1$ is chosen such that $\bbP\left(\bsTheta^\mathrm{OPT}\left(\hat{\bfx}, \bfu\right) \cap \mC(\bar{\bstheta}, \alpha_1) \right) = 1$, we have
\begin{align*}
    \mathrm{POG}(\bar{\bfx}_\mathrm{CIO}) 
    = \ 
    & \bbE
        \left[
            f\left(\hat{\bstheta}, \bar{\bfx}_\mathrm{CIO}(\bfu)\right) 
            - f\left(\hat{\bstheta}, \hat{\bfx} \right) 
        \right] \\
    \leq \ 
    & \bbE
        \left\{
            \nu(\hat{\bfx})(\eta - 2\cos2\alpha_1 + 2)
            + \eta \nu\left[\bar{\bfx}_\mathrm{CIO}(\bfu)\right]
        \right\} \\
    = \ 
    & \mu(\eta - 2\cos2\alpha_1 + 2)
        + \eta \mu_\mathrm{CIO}
\end{align*}
where $\mu := \bbE\left[\nu(\hat{\bfx}) \right]$ and $\mu_\mathrm{CIO} := \bbE \left(\nu\left[ \bar{\bfx}_\textrm{CIO}(\bfu) \right]\right)$.
\endproof

\begin{proposition}[Conformal IO Achieves Bounded AOG] \label{prop:aog}
   Let $\bar{\bfx}_\mathrm{CIO}(\bfu)$ be an optimal solution to $\mathbf{RFO}\left(\mC(\bar{\bstheta}, \alpha_1), \bfu \right)$ for any $\bfu\in \mU$, where $\bar{\bstheta} \in \bbR^d$ and $\alpha_1$ are chosen such that, for a new sample $(\bstheta', \bfu')$ from $\bbP_{(\bstheta, \bfu)}$ and $\bfx' = \tilde{\bfx}(\bstheta', \bfu')$, $\bbP \left( \mC(\bar{\bstheta}, \alpha_1) \cap \bsTheta^\mathrm{OPT}(\bfu', \bfx') \neq \emptyset \right) = 1$. If Assumptions \ref{asm:bounded_inv_feas}--\ref{asm:bnd_div} hold, then
    \begin{equation}
        \mathrm{AOG}(\bar{\bfx}_\mathrm{CIO}) 
        \leq 
        (2 - 2 \cos 2\alpha_1 + \eta + \sigma) \mu^*   
        + (\eta + \sigma) \mu_\mathrm{CIO}
    \end{equation}
    where $\mu^* := \bbE\left(\nu[ \tilde{\bfx}(\bstheta^*, \bfu)] \right)$.
\end{proposition}

\proof{} 
We first derive an upper bound on the optimality gap of the suggested decision $\bar{\bfx}_\textrm{CIO}(\bfu)$ as evaluated using $\bstheta^*$ for any $\bfu\in \mU$. Let $(\hat{\bstheta}, \bfu)$ be a sample from $\bbP_{(\bstheta, \bfu)}$, $\hat{\bfx} = \tilde{\bfx}(\hat{\bstheta}, \bfu)$, and $\tilde{\bstheta}$ be an element of $\bsTheta^\mathrm{OPT}\left(\hat{\bfx}, \bfu\right) \cap \mC(\bar{\bstheta}, \alpha_1)$, which is non-empty almost surely because $\alpha_1$ is chosen such that $\bbP\left(\bsTheta^\mathrm{OPT}\left(\hat{\bfx}, \bfu\right) \cap \mC(\bar{\bstheta}, \alpha_1) \right) = 1$. 
Let $\bar{\bstheta}_\textrm{CIO}(\bfu)$ denote the optimal solution to the inner maximization problem in $\mathbf{RFO}\left(\mC(\bar{\bstheta}, \alpha_1), \bfu\right)$ when the outer decision variable is set to $\bar{\bfx}_\mathrm{CIO}(\bfu)$. For any $\bfu\in \mU$, let $\bfx^*(\bfu) := \tilde{\bfx}(\bstheta^*, \bfu)$
and 
$\bstheta^*_\textrm{CIO}(\bfu)$ denote the optimal solution to the inner maximization problem in $\mathbf{RFO}\left(\mC(\bar{\bstheta}, \alpha_1), \bfu\right)$ when the outer decision variable is set to $\bfx^* (\bfu)$,
we have

\begin{align*}
    & f\left(\bstheta^*, \bar{\bfx}_\textrm{CIO}(\bfu)\right)
    - f\left(\bstheta^*, \bfx^*(\bfu) \right) \\
    \leq \ 
    & f\left(\bbE(\hat{\bstheta}), \bar{\bfx}_\textrm{CIO}(\bfu)\right)
        - f\left(\bbE(\hat{\bstheta}), \bfx^*(\bfu) \right) 
        + \left(
            \nu\left[ \bar{\bfx}_\textrm{CIO}(\bfu) \right]
            + \nu\left[ \bfx^*(\bfu) \right]
            \right)
            \|\bstheta^* -  \bbE(\hat{\bstheta})\|_2  \\
    \leq \ 
    & f\left(\bbE(\hat{\bstheta}), \bar{\bfx}_\textrm{CIO}(\bfu)\right)
        - f\left(\bbE(\hat{\bstheta}), \bfx^*(\bfu) \right) 
        + \sigma  
            \left(
            \nu\left[ \bar{\bfx}_\textrm{CIO}(\bfu) \right]
            + \nu\left[ \bfx^*(\bfu) \right]
            \right) \\
    = \ 
    & \bbE\left[
        f\left(\hat{\bstheta}, \bar{\bfx}_\textrm{CIO}(\bfu)\right)
        - f\left(\hat{\bstheta}, \bfx^*(\bfu) \right) 
        \right]
        + \sigma  
            \left(
            \nu\left[ \bar{\bfx}_\textrm{CIO}(\bfu) \right]
            + \nu\left( \bfx^*(\bfu) \right)
            \right) \\
    \leq \ 
    & \bbE\left[
        f \left(\tilde{\bstheta}, \bar{\bfx}_\textrm{CIO}(\bfu)\right)
        - f\left(\tilde{\bstheta}, \bfx^*(\bfu) \right)
        + (\nu\left[\bar{\bfx}_\mathrm{CIO}(\bfu)\right] 
            + \nu\left[\hat{\bfx}\right])
        \|\hat{\bstheta} - \tilde{\bstheta} \|_2
        \right] \\
    & + \sigma  
            \left(
            \nu\left[ \bar{\bfx}_\textrm{CIO}(\bfu) \right]
            + \nu\left( \bfx^*(\bfu) \right)
            \right) \\
    \leq \ 
    & \bbE\left[
        f \left(\tilde{\bstheta}, \bar{\bfx}_\textrm{CIO}(\bfu)\right)
        - f\left(\tilde{\bstheta}, \bfx^*(\bfu) \right)
        + (\nu\left[\bar{\bfx}_\mathrm{CIO}(\bfu)\right] 
            + \nu\left[\hat{\bfx}\right])
        \eta
        \right]\\
    & + \sigma  
            \left(
            \nu\left[ \bar{\bfx}_\textrm{CIO}(\bfu) \right]
            + \nu\left( \bfx^*(\bfu) \right)
            \right) \\
    \leq \ 
    & \bbE\left[
        f \left(\bar{\bstheta}_\textrm{CIO}(\bfu), \bar{\bfx}_\textrm{CIO}(\bfu)\right)
        - f\left(\tilde{\bstheta}, \bfx^*(\bfu) \right)
        \right]
        + (\eta + \sigma) \left( 
            \nu\left[ \bar{\bfx}_\textrm{CIO}(\bfu)\right] 
            + \nu\left( \bfx^*(\bfu) \right)
            \right) \\
    \leq \ 
    & \bbE\left[
        f \left(\bstheta^*_\textrm{CIO}(\bfu), \bfx^*(\bfu) \right)
        - f\left(\tilde{\bstheta}, \bfx^*(\bfu) \right)
        \right]
        + (\eta + \sigma) \left( 
            \nu\left[ \bar{\bfx}_\textrm{CIO}(\bfu)\right] 
            + \nu\left( \bfx^*(\bfu) \right)
            \right) \\
    \leq \ 
    & \bbE\left[
        \nu(\bfx^*(\bfu)) \|\bstheta^*_\textrm{CIO}(\bfu) - \tilde{\bstheta} \|_2
        \right]
        + (\eta + \sigma) \left( 
            \nu\left[ \bar{\bfx}_\textrm{CIO}(\bfu)\right] 
            + \nu\left( \bfx^*(\bfu) \right)
            \right) \\
    \leq \ 
    & 2 \nu\left(\bfx^*(\bfu)\right) (1 - \cos 2\alpha_1)
        + (\eta + \sigma) \left( 
            \nu\left[ \bar{\bfx}_\textrm{CIO}(\bfu)\right] 
            + \nu\left( \bfx^*(\bfu) \right)
            \right) \\
    \leq \ 
    & (2 - 2 \cos 2\alpha_1 + \eta + \sigma)             
        \nu\left(\bfx^*(\bfu)\right)
        + (\eta + \sigma)
            \nu\left[ \bar{\bfx}_\textrm{CIO}(\bfu) \right]
\end{align*}

The first line holds because of Lemma \ref{lemma:lipschiz}.
The second line holds due to Assumption \ref{asm:bnd_div}. 
The third line holds because $f$ is linear in $\bstheta$. The expectation is taken over $\bbP_{\bstheta}$. 
The fourth line holds due to Lemma \ref{lemma:lipschiz}. 
The fifth line holds due to Assumption \ref{asm:bounded_inv_feas}. 
The sixth line holds because of the definition of $\bar{\bstheta}_\textrm{CIO}(\bfu)$. 
The seventh line holds because $\left( \bar{\bfx}_\mathrm{CIO}(\bfu), \bar{\bstheta}_\textrm{CIO}(\bfu) \right)$ is an optimal solution to $\mathbf{RFO}\left(\mC(\bar{\bstheta}, \alpha_1), \bfu\right)$.
The eigth line holds due to Lemma \ref{lemma:lipschiz}. 
The ninth line holds since both $\bstheta^*_\textrm{CIO}(\bfu)$ and $\tilde{\bstheta}$ are on the unit sphere and the angle between them is no greater than $2\alpha_1$, then the $L_2$ distance between them is upper bounded by $2(1 - \cos2\alpha_1)$.

Next, we bound the $\mathrm{AOG}$ of $\bar{\bfx}_\textrm{CIO}$. We have
\begin{align*}
    \mathrm{AOG}(\bar{\bfx}_\textrm{CIO})
    =\ 
    & \bbE
    \left[
        f\left(\bstheta^*, \bar{\bfx}_\textrm{CIO}(\bfu)\right)
        - f\left(\bstheta^*, \bfx^*(\bfu) \right)
    \right] \\
    \leq \ 
    & \bbE
    \left[
        (2 - 2 \cos 2\alpha_1 + \eta + \sigma)             
        \nu\left(\bfx^*(\bfu)\right)
        + (\eta + \sigma)
            \nu\left[ \bar{\bfx}_\textrm{CIO}(\bfu) \right]
    \right] \\
    = \ 
    & (2 - 2 \cos 2\alpha_1 + \eta + \sigma) \mu^*   
        + (\eta + \sigma) \mu_\mathrm{CIO}
\end{align*}
where $\mu^* := \bbE\left(\nu[ \bfx^*(\bfu) ]\right)$.
\endproof



\section{Numerical Experiment Details} \label{ec:implementation}

\subsection{Computational Setup} \label{subsec:comp_setup}
All the algorithms are implemented and test using Python 3.9.1 on a MacBook Pro with an Apple M1 Pro processor and 16 GB of RAM. Optimization models are implemented with Gurobi 9.5.2.

\subsection{Forward Problems} \label{subsec:fwd_prob}

\subsubsection{Shortest-path}

We consider the shortest path problem on a 5$\times$5 grid network $G(\mN, \mE)$ where $\mN$ and $\mE$ indicate the node and edge sets, respectively. Let $\mE^+(i)$ and $\mE^-(i)$ denote the sets of edges that enter and leave node $i\in \mN$, respectively. Let $u^o$ and $u^d$ denote the origin and destination of the trip, respectively. We define $x_{ij} \in \mE$ as binary decision variables that take 1 if road $(i, j)$ is traversed for any $(i, j)\in \mE$. The shortest path problem is presented as follows.
\begin{subequations}
\begin{align}
    \underset{\bfx}{\textrm{minimize}} \quad
        & \sum_{(i, j)\in \mE} \theta_{ij} x_{ij} \\
    \textrm{subject to} \quad
        & \sum_{(j, i)\in \mE^+(i)} x_{ji}
            - \sum_{(i, j)\in \mE^-(i)} x_{ij}
            = 
            \begin{cases}
                1,  & \textrm{ if } i = u^d \\
                -1, & \textrm{ if } i = o^d\\
                0,  & \textrm{ otherwise}
            \end{cases}
            , \quad
            \forall i \in \mN\\
        & x_{ij} \in \{0, 1\},
            \quad (i, j)\in \mE.
\end{align}
\end{subequations}
The objective function minimizes the total travel cost. The first set of constraints are the flow-balancing constraints that make sure we can find a path from $u_o$ to $u_d$. The second set of constraints specify the range of our decision variables. Note that the constriant matrix is totally unimodular, so we can replace the binary constraints with $0\leq x_{ij} \leq 1$ for any $(i, j)\in \mE$ when implementing this model.

\subsubsection{Knapsack}
We consider a knapsack problem of $d=10$ items. We define binary decision variables $x_i$ that indicate if item $i\in [d]$ is selected ($=1$) or not ($=0$). The knapsack problem is presented as follows.
\begin{subequations}
\begin{align}
    \underset{\bfx}{\textrm{maximize}}\quad 
        & \sum_{i\in [d]} \theta_i x_i \\
    \textrm{subject to} \quad
        & \sum_{i\in [d]} w_i x_i \leq u\\
        & x_i\in \{0, 1\}, \forall i \in [d].
\end{align}
\end{subequations}

The objective maximizes the total value of the selected items. The first constraint enforces a total budget for item selection. The second set of constraints specify the range of our decision variables.

\subsection{Obtaining a Point Estimation} \label{ec:cp_ddio}

We consider two methods to obtain point estimations of the unknown parameters. They are i) data-driven inverse optimization with the sub-optimality loss and ii) the gradient-based method proposed by \citet{berthet2020learn}. We implement the method from \citet{berthet2020learn} with the package provided by \citet{tang2022pyepo}. Hyper-parameters are tuned using a separate validation set of 200 decision data points. Batch size is set to 64. We use the Adam optimizer with an initial learning rate of 0.1. We train the model for 20 epochs.

We present the implementation details of the data-driven inverse optimization method next. We consider solving
\begin{subequations}
\label{ddio:master}
\begin{align}
    \underset{\bstheta \in \bbR^{|\mE|}, \bsepsilon \in \bbR_+^{n_\textrm{train}}}{\textrm{minimize}} \quad
        & \frac{1}{N_\textrm{train}}
        \sum_{k\in \mK_\textrm{train}}
            l_k \\
    \textrm{subject to} \quad
        \label{ddio:opt_gap}
        & l_k \geq \bstheta^\intercal \hat{\bfx}_k
            - \bstheta^\intercal \bfx,
            \quad \forall \bfx \in \mX_k,\ 
            k\in \mK_\textrm{train} \\
        & \| \bstheta - \mathbf{1}\|_1 \leq \frac{|\mE|}{4}.
\end{align}
\end{subequations}
This problem is initialized without Constraints~\eqref{ddio:opt_gap}, which were added iteratively using a cutting-plane method. Specifically, in each iteration, after solving Problem~\eqref{ddio:master}, let $\bstheta'$ and $\{l'_k\}_{k\in \mK_\textrm{train}}$ be the optimal solution. For each data point $k\in \mK_\textrm{train}$, we solve the following sub-problem
\begin{align}
    \underset{\bfx_k \in \mX(\bfu_k)}{\textrm{minimize}} \quad 
        \bstheta'^\intercal \bfx_k.
\end{align}
Let $\bfx_k'$ be the optimal solution to the sub-problem. If $l_k' < \bstheta'^\intercal \hat{\bfx}_k - \bstheta'^\intercal \bfx_k'$, we add the following cut to Problem~\eqref{ddio:master}
\begin{equation}
    l_k 
    \geq 
    \bstheta^\intercal \hat{\bfx}_k
        - \bstheta^\intercal \bfx'_k.
\end{equation}

We keep running this procedure until no cut is added to the master Problem~\eqref{ddio:master}.

\subsection{Solving the Calibration Problem} \label{ec:calib}

\subsubsection{Shortest Path} \label{ec:calib_sp}

For each data point in the validation set, we calculate the value of $c_k$ by solving the following problem
\begin{subequations}
\begin{align}
    \underset{\bstheta \in \bbR^{|\mE|}, \bfw \in \bbR^{\mN}, \bfv \in \bbR_+^{|\mE|}}{\textrm{maximize}}\quad 
        & \bar{\bstheta}^\intercal \bstheta \\
    \textrm{subject to} \quad
        & w_{d_k} - w_{o_k} - \sum_{(i, j)\in \mE} v_{ij}
            = \bstheta^\intercal \hat{\bfx}_k \\
        & w_j - w_{i} - v_{ij} \leq c_{ij},
            \quad \forall (i, j)\in \mE \\
        & \|\bstheta\|_2 \leq 1.
\end{align}
\end{subequations}
where $\bfw\in \bbR^\mN$ and $\bfv\in \bbR_+^{|\mE|}$, respectively, denote the dual variables associated with the flow-balancing constraints and the capacity constraints in the primal problem. The first constraint enforces strong duality. The second set of constraints are the dual feasibility constraints. The last constraint ensures the optimal solution is on the unit sphere. Note that we do not need to enforce $\|\bstheta\|_2 = 1$ because this is a maximization problem.

\subsubsection{Knapsack} \label{ec:calib_knapsack}

For each data point in the validation set, we calculate the value of $c_k$ by solving the following calibration problem
\begin{subequations}
\begin{align}
    \underset{\bstheta \in \bbR^{d}}{\textrm{maximize}}\quad 
        & \bar{\bstheta}^\intercal \bstheta \\
    \textrm{subject to} \quad
        \label{ks_con:opt_con}
        & \bstheta^\intercal \hat{\bfx}_k \geq \bstheta^\intercal \bfx,
            \quad \forall \bfx\in \mX(\bfu_k) \\
        & \|\bstheta\|_2 \leq 1.
\end{align}
\end{subequations}
We initialize this problem without Constraints~\eqref{ks_con:opt_con}. In each iteration, after solving the calibration problem, let $\bstheta'$ denote the optimal solution. We solve $\mathbf{FO}(\bstheta', \bfu_k)$ and let $\bfx'$ denote the optimal solution. If $\bstheta'^\intercal\bfx' > \bstheta'^\intercal \hat{\bfx}_k$, we then add the corresponding cut to the model. We keep running this process until no cut is added.

\subsection{Solving the Robust Forward Problem} \label{ec:robust_fwd}

Let $\alpha = \cos^{-1}\left(\Gamma_k(\{c_k\}_{k\in \mK_\textrm{val}}) \right)$. We next solve the following robust model to recommend a new decision to prescribe a decision given a $u\in \mU$.
\begin{subequations}
\begin{align}
    \underset{\bfx \in \mX(\bfu)}{\textrm{minimize}}\ 
    \underset{\bstheta \in \bbR^{|\mE|}}{\textrm{maximize}} \quad 
        & \bstheta^\intercal \bfx \\
    \textrm{subject to} \quad
        & \bar{\bstheta}^\intercal \bstheta \geq \cos(\alpha) \\
        & \|\bstheta\|_2 \leq 1.
\end{align}
\end{subequations}

We initialize this problem as follows.
\begin{subequations} \label{rfo:master}
\begin{align}
    \underset{\bfx \in \mX(\bfu), \Omega\in \bbR_+}{\textrm{minimize}}\ 
    \Omega \\
    \textrm{subject to} \quad
        & \bstheta^\intercal \bfx \leq \Omega,
            \quad \forall \bstheta \in \tilde{\bsTheta}.
\end{align}
\end{subequations}

We initialize $\tilde{\bsTheta} = \emptyset$. We first solve Problem~\eqref{rfo:master}, let $\bfx'$ and $\Omega'$ denote the optimal solution. Then we solve the following sub-problem
\begin{subequations} \label{rfo:sub}
\begin{align}
     \underset{\bstheta \in \bbR^{|\mE|}}{\textrm{maximize}} \quad 
        & \bstheta^\intercal \bfx' \\
    \textrm{subject to} \quad
        & \bar{\bstheta}^\intercal \bstheta \geq \cos(\alpha) \\
        & \|\bstheta\|_2 \leq 1.
\end{align}
\end{subequations}

Let $\bstheta'$ denote the optimal solution to the sub-problem. If $\bstheta^\intercal \bfx' > \Omega'$, then we add $\bstheta'$ to $\tilde{\bsTheta}$ and re-solve Problem~\eqref{rfo:master}. We keep running this procedure until no new solution is added to $\tilde{\bsTheta}$.

\end{document}